# Beta generated Kumaraswamy-*G* and other new families of distributions

Laba Handique and Subrata Chakraborty*

Department of Statistics, Dibrugarh University

Dibrugarh-786004, India

*Corresponding Author. Email: subrata_stats@dibru.ac.in

(August 21, 2016)

**Abstract**

A new generalization of the family of Kumaraswamy-*G* (Cordeiro and de Castro, 2011) distribution that includes three recently proposed families namely the Garhy generated family (Elgarhy *et al.,* 2016), Beta-Dagum and Beta-Singh-Maddala distribution (Domma and Condino, 2016) is proposed by constructing beta generated Kumaraswamy-*G* distribution. Useful expansions of the pdf and the cdf of the proposed family is derived and seen as infinite mixtures of the Kumaraswamy-*G* distribution. Order statistics, Probability weighted moments, moment generating function, Rényi entropies, quantile power series, random sample generation, asymptotes and shapes are also investigated. Two methods of parameter estimation are presented. Suitability of the proposed model in comparisons to its sub models is carried out considering two real life data sets. Finally, some new classes of beta generated families are proposed for future investigations.

**Key words**: *Beta Generated family; Kumaraswamy-G distribution; Exponentiated family; Beta Marshall-Olkin Kumaraswamy-G; Beta Kumaraswamy Marshall-Olkin-G; Beta Generalized Marshall-Olkin Kumaraswamy-G and Beta Kumaraswamy Generalized Marshall-Olkin-G.*

# 1. Introduction

Here we briefly introduce the Beta-$G$ (Eugene et al., 2002 and Jones, 2004) and Kumaraswamy-$G$ ($Kw-G$) (Cordeiro and de Castro, 2011) family of distributions.

## 1.1 Beta-$G$ family of distributions

For a given distribution with pdf $f(t)$ and cdf $F(t)$, the cdf of beta-$G$ (Eugene et al., 2002 and Jones, 2004) family of distribution is given respectively by

$$F^{BG}(t;m,n) = \frac{1}{B(m,n)} \int_o^{F(t)} v^{m-1}(1-v)^{n-1} dv = \frac{B_{F(t)}(m,n)}{B(m,n)} = I_{F(t)}(m,n) \qquad (1)$$

and

$$f^{BG}(t;m,n) = \frac{1}{B(m,n)} f(t) F(t)^{m-1}[1-F(t)]^{n-1} = \frac{1}{B(m,n)} f(t) F(t)^{m-1} \overline{F}(t)^{n-1} \qquad (2)$$

Where $I_t(m,n) = B(m,n)^{-1} \int_0^t x^{m-1}(1-x)^{n-1} dx$ denotes the incomplete beta function ratio. The survival function (sf), hazard rate function (hrf), reverse hazard rate function (rhrf) and cumulative hazard rate function (chrf) are given respectively by

$$\overline{F}^{BG}(t;m,n) = P[T>t] = 1 - I_{F(t)}(m,n) = \frac{B(m,n) - B_{F(t)}(m,n)}{B(m,n)}$$

$$h^{BG}(t;m,n) = f^{BG}(t;m,n) / \overline{F}^{BG}(t;m,n) = \frac{f(t) F(t)^{m-1} \overline{F}(t)^{n-1}}{B(m,n) - B_{F(t)}(m,n)}$$

$$r^{BG}(t;m,n) = f^{BG}(t;m,n) / F^{BG}(t;m,n) = \frac{f(t) F(t)^{m-1} \overline{F}(t)^{n-1}}{B_{F(t)}(m,n)}$$

and $H^{BG}(t;m,n) = -\log[\overline{F}(t;m,n)] = -\log\left[\frac{B(m,n) - B_{F(t)}(m,n)}{B(m,n)}\right]$.

Some of the well known works related to beta-generated (beta-$G$) family (Jones 2004), are the beta Gumbel (Nadarajah and Kotz, 2004) beta Frechet distribution (Nadarajah and Gupta, 2004), beta generalized exponential (Barreto-Souza *et al.*, 2010) the beta extended-$G$ family (Cordeiro *et al.*, 2012), Kumaraswamy beta generalized family (Pescim *et al.*, 2012), beta generalized Weibull distribution (Singla *et al.*, 2012), beta generalized Rayleigh distribution (Cordeiro *et al.*, 2013), beta extended half normal distribution (Cordeiro *et al.*, 2014), beta log-logistic distribution (Lemonte, 2014), beta generalized inverse Weibull distribution (Baharith *et al.*, 2014), beta Marshall-Olkin family of distribution (Alizadeh et al., 2015), beta exponential Frechet distribution (Mead *et al.*, 2016), beta-Dagum distribution and beta-Singh-Maddala distribution (Domma and Condino, 2016) among others.

## 1.2 Kumaraswamy $G$ ($Kw-G$) family of distributions

For a baseline cdf $G(t)$ with pdf $g(t)$, Cordeiro and de Castro (2011) defined $Kw-G$ distribution with respective cdf and pdf

$$F^{KwG}(t;a,b) = 1 - [1 - G(t)^a]^b \tag{3}$$

and

$$f^{KwG}(t;a,b) = a\,b\,g(t)G(t)^{a-1}[1-G(t)^a]^{b-1} \tag{4}$$

Where $t > 0$, $g(t) = \frac{d}{dt}G(t)$ and $a > 0, b > 0$ are shape parameters in addition to those in the baseline distribution. Corresponding sf, hrf, rhrf and chrf are respectively given by

$$\bar{F}^{KwG}(t;a,b) = 1 - F^{KwG}(t) = [1-G(t)^a]^b \tag{5}$$

$$h^{KwG}(t;a,b) = f^{KwG}(t;a,b)/\bar{F}^{KwG}(t;a,b)$$

$$= a\,b\,g(t)G(t)^{a-1}[1-G(t)^a]^{b-1}/[1-G(t)^a]^b$$

$$= a\,b\,g(t)G(t)^{a-1}[1-G(t)^a]^{-1}$$

$$r^{KwG}(t;a,b) = f^{KwG}(t;a,b)/F^{KwG}(t;a,b)$$

$$= a\,b\,g(t)G(t)^{a-1}[1-G(t)^a]^{b-1}/1-[1-G(t)^a]^b$$

$$= a\,b\,g(t)G(t)^{a-1}[1-G(t)^a]^{b-1}\{1-[1-G(t)^a]^b\}^{-1}$$

and $H^{KwG}(t;a,b) = -b\log[1-G(t)^a]$ respectively.

Some of the notable distributions derived under the scheme of $Kw-G$ family (Cordeiro and de Castro, 2011), Kumaraswamy Normal distribution (Correa *et al.*, 2012), Kumaraswamy generalized Pareto distribution (Nadarajah and Eljbri, 2013), Kumaraswamy linear exponential distribution (Elbatal, 2013), Kumarswamy Lomax distribution (Shams, 2013), Kumaraswamy modified Weibull distribution (Cordeiro *et al.*, 2014), Kumaraswamy generalized Rayleigh distribution (Gomes et al., 2014), Marshall-Olkin Kumaraswamy-*G* (Handique and Chakraborty, 2015a) Generalized Marshall-Olkin Kumaraswamy-*G* family of distribution (Handique and Chakraborty, 2015b), Kumaraswamy Marshall-Olkin family of distribution (Alizadeh *et al.*, 2015) and Kumaraswamy Pareto IV Distribution (Tahir *et al.*, 2015) among others.

In this article we propose a family of beta generated $Kw-G$ distribution by considering the $Kw-G$ (Cordeiro and de Castro, 2011) family as the base line distribution in the beta-$G$ (Eugene et al., 2002 and Jones, 2004) family and investigate some of its general properties. The rest of this article is organized in eight sections. In section 2 the new family is defined. Important special cases of the family along with their shape and main reliability characteristics are presented in the next section. In section 4 we discuss some general results of the family. Different methods of estimation of parameters are presented in section 5. Application of the proposed family is considered in section 6. In section 7 we propose a few more families of beta generated distributions for future studies. The paper ends with a conclusion in section 8 followed by an appendix to derive asymptotic confidence bounds.

**2. A new generalization: Beta Kumarswamy-$G$ ($BKw-G$) family of distributions**

Here we propose a new beta generated family by considering the cdf, pdf and sf of $Kw-G$ (Cordeiro and de Castro, 2011) distribution in (3), (4) and (5) as the $F(t)$, $f(t)$ and $\bar{F}(t)$ respectively in the beta-$G$ (Eugene et al., 2002 and Jones, 2004) formulation in (2) and call it $BKw-G$ distribution. The pdf and cdf of $BKw-G$ are given respectively by

$$f^{BKwG}(t;a,b,m,n) = \frac{ab\,g(t)G(t)^{a-1}[1-G(t)^a]^{bn-1}[1-[1-G(t)^a]^b]^{m-1}}{B(m,n)}, \qquad (6)$$

$$0 < t < \infty,\; 0 < a,b < \infty,\; m,n > 0$$

and $\qquad F^{BKwG}(t;a,b,m,n) = I_{1-[1-G(t)^a]^b}(m,n) \qquad (7)$

The sf, hrf, rhrf and chrf of $BKw-G$ distribution are respectively obtained as

$$\bar{F}^{BKwG}(t;a,b,m,n) = 1 - I_{1-[1-G(t)^a]^b}(m,n)$$

$$h^{BKwG}(t;a,b,m,n) = \frac{ab\,g(t)G(t)^{a-1}[1-G(t)^a]^{bn-1}[1-[1-G(t)^a]^b]^{m-1}}{B(m,n)[1-I_{1-[1-G(t)^a]^b}(m,n)]} \qquad (8)$$

$$r^{BKwG}(t;a,b,m,n) = \frac{ab\,g(t)G(t)^{a-1}[1-G(t)^a]^{bn-1}[1-[1-G(t)^a]^b]^{m-1}}{B(m,n)\,I_{1-[1-G(t)^a]^b}(m,n)} \qquad (9)$$

$$H^{BKwG}(t;a,b,m,n) = -\log[1 - I_{1-[1-G(t)^a]^b}(m,n)]$$

This new family can also be referred to as the new generalized $Kw-G$ ($GKw-G$) family of distribution. For $n=1$, it reduces to the recently proposed Garhy generated family (Garhy et al., 2016) that is $BKw\text{-}G(m,1,a,b) \equiv GH\text{-}G(m,a,b)$ with cdf and pdf respectively are given by $F^{GHG}(t;\alpha,\beta,\gamma) = [1-[1-G(t)^\alpha]^\beta]^\gamma$ and

$f^{GHG}(t;\alpha,\beta,\gamma) = \alpha\beta\gamma\,g(t)G(t)^{\alpha-1}[1-G(t)^\alpha]^{\beta-1}[1-[1-G(t)^\alpha]^\beta]^{\gamma-1}\; t>0,$ Where $\alpha,\beta,\gamma > 0$ are three shape parameters.

For $m=n=1$, $BKw\text{-}G(1,1,a,b) \equiv Kw\text{-}G(a,b)$, (Cordeiro and de Castro, 2011) and if $a=b=1$ then $BKw\text{-}G(m,n,1,1) \equiv B(m,n)$, (Eugene et al., 2002 and Jones, 2004).

**2.1 Genesis of the distribution**

If $m$ and $n$ are both integers, then the probability distribution of $BKw-G$ arises as distribution of the $m^{th}$ order statistics from a random sample of size $m+n-1$ from distribution.

**Proof:** Let $T_1, T_2, ..., T_{m+n-1}$ be a random sample from a population with cdf $1 - [1 - G(t)^a]^b$. Then the pdf of the $m^{th}$ order statistics $T_{(m)}$ is given by

$$= \frac{(m+n-1)!}{(m-1)! \, [(m+n-1)-m]!} [1 - [1-G(t)^a]^b]^{m-1} [\{1-G(t)^a\}^b]^{(m+n-1)-m}$$

$$\times a b \, g(t) G(t)^{a-1} [1-G(t)^a]^{b-1}$$

$$= \frac{\Gamma(m+n)}{\Gamma(m) \, \Gamma(n)} [1 - [1-G(t)^a]^b]^{m-1} [\{1-G(t)^a\}^b]^{n-1}$$

$$\times a b \, g(t) G(t)^{a-1} [1-G(t)^a]^{b-1}$$

$$= \frac{a b \, g(t)}{B(m,n)} G(t)^{a-1} [1-G(t)^a]^{bn-1} [1 - [1-G(t)^a]^b]^{m-1}$$

## 2.2 Shape of the density and hazard functions

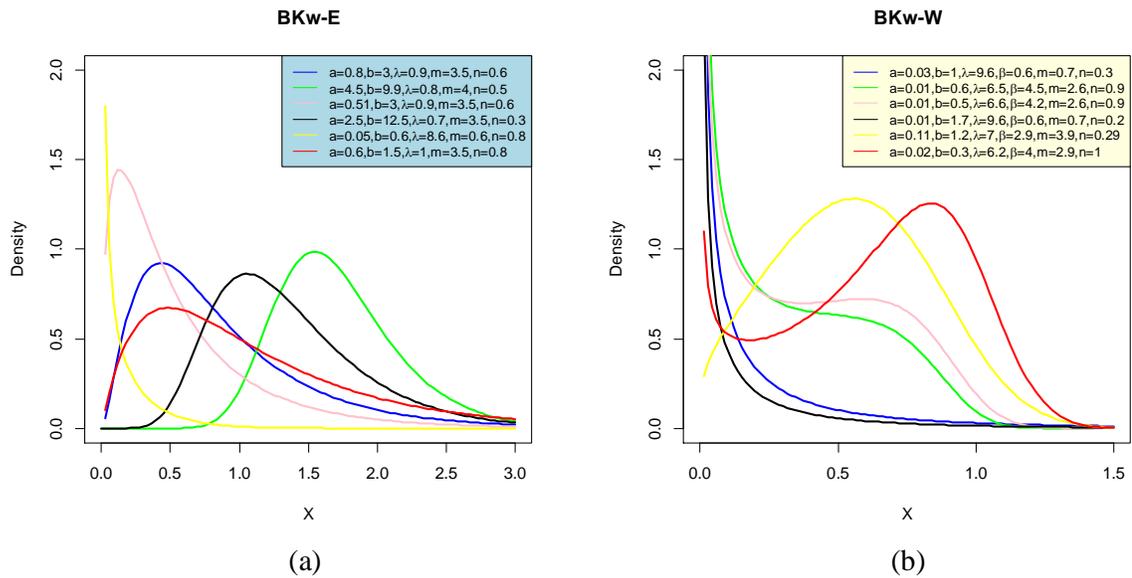

(a)   (b)

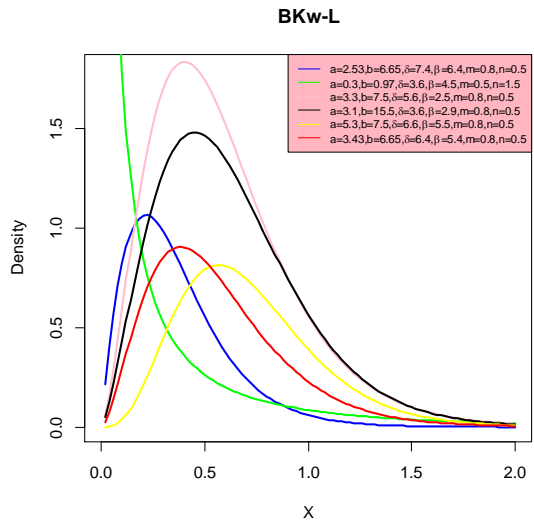
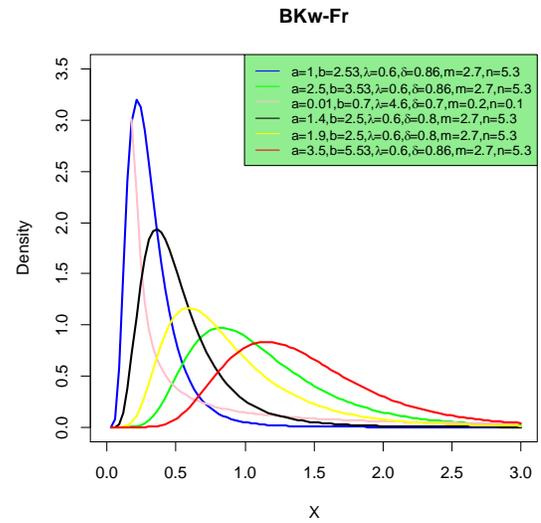

(c)                                        (d)

**Fig 1:** Density plots of (a) $BKw-E$, (b) $BKw-W$, (c) $BKw-L$ and (d) $BKw-Fr$ distributions.

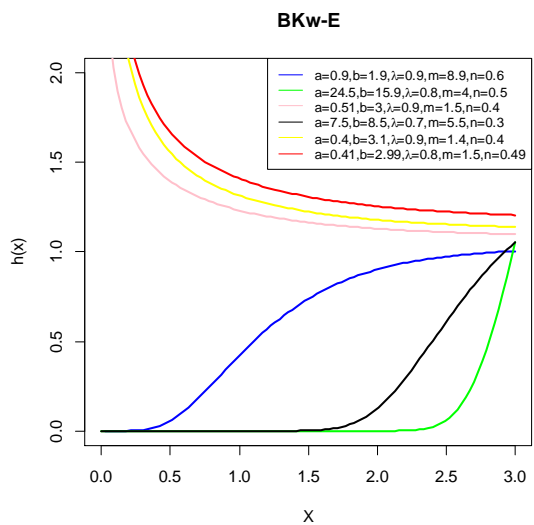
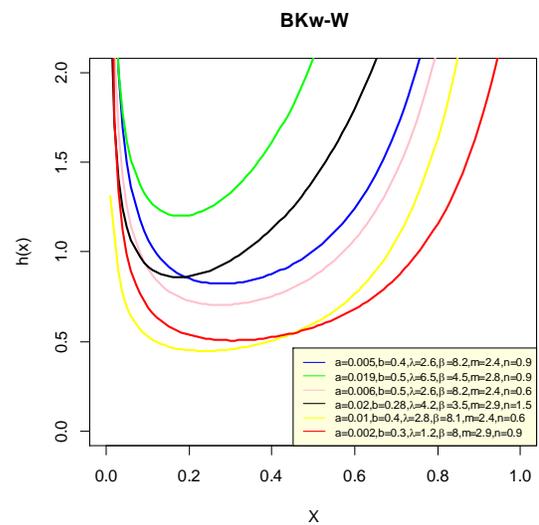

(a)                                        (b)

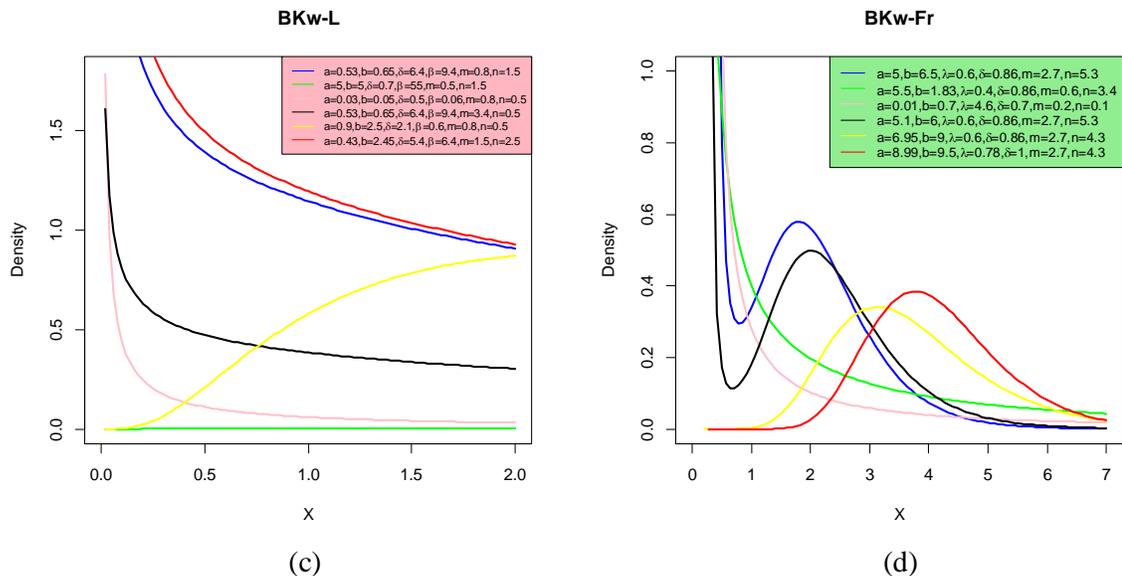

(c)                                                   (d)

**Fig 2:** Hazard plots of (a) $BKw-E$, (b) $BKw-W$, (c) $BKw-L$ and (d) $BKw-Fr$ distributions.

From the plots in figure 1 and 2 it can be seen that the family is very flexible and can offer many different types of shapes. It offers IFR, DFR even bath tub shaped hazard rate.

### 3. Some special $BKw-G$ distributions

A number of known and new distributions can be obtained as special cases of the $BKw-G$ family of distributions. Some of which are listed next.

3.1 The $BKw$- Dagum distribution ($BKw-Da$) distribution:

Suppose the base line distribution is the Dagum distribution (Dagum, 1977, 1980) with cdf and pdf given by

$$G(t;\tau) = (1+\tau_2 t^{-\tau_3})^{-\tau_1} \quad \text{and} \quad g(t;\tau) = \tau_1 \tau_2 \tau_3 t^{-(\tau_3+1)}(1+\tau_2 t^{-\tau_3})^{-\tau_1-1}$$

with $\tau' = (\tau_1, \tau_2, \tau_3)$ and $\tau_i > 0$ for $i=1,2,3$ where $\tau_2$ is a scale parameter and $\tau_1$ and $\tau_3$ are shape parameters respectively, then the corresponding pdf and cdf of $BKw-Da$ distribution becomes

$$f^{BKwDa}(t;a,b,m,n,\tau) = \frac{\left[\begin{array}{l} ab\tau_1\tau_2\tau_3 t^{-(\tau_3+1)}(1+\tau_2 t^{-\tau_3})^{-\tau_1-1}\{(1+\tau_2 t^{-\tau_3})^{-\tau_1}\}^{a-1} \\ [1-\{(1+\tau_2 t^{-\tau_3})^{-\tau_1}\}^a]^{bn-1}[1-[1-\{(1+\tau_2 t^{-\tau_3})^{-\tau_1}\}^a]^b]^{m-1} \end{array}\right]}{B(m,n)}$$

and $\quad F^{BKwDa}(t;a,b,m,n,\tau) = I_{1-[1-\{(1+\tau_2 t^{-\tau_3})^{-\tau_1}\}^a]^b}(m,n).$

sf: $\quad \overline{F}^{BKwDa}(t;a,b,m,n,\tau) = 1 - I_{1-[1-\{(1+\tau_2 t^{-\tau_3})^{-\tau_1}\}^a]^b}(m,n)$

hrf: $h^{BKwDa}(t;a,b,m,n,\tau) = \dfrac{\begin{bmatrix} ab\tau_1\tau_2\tau_3 t^{-(\tau_3+1)}(1+\tau_2 t^{-\tau_3})^{-\tau_1-1}\{(1+\tau_2 t^{-\tau_3})^{-\tau_1}\}^{a-1} \\ [1-\{(1+\tau_2 t^{-\tau_3})^{-\tau_1}\}^a]^{bn-1}[1-[1-\{(1+\tau_2 t^{-\tau_3})^{-\tau_1}\}^a]^b]^{m-1} \end{bmatrix}}{B(m,n)\,[1-I_{1-[1-\{(1+\tau_2 t^{-\tau_3})^{-\tau_1}\}^a]^b}(m,n)]}$

rhrf: $r^{BKwDa}(t;a,b,m,n,\tau) = \dfrac{\begin{bmatrix} ab\tau_1\tau_2\tau_3 t^{-(\tau_3+1)}(1+\tau_2 t^{-\tau_3})^{-\tau_1-1}\{(1+\tau_2 t^{-\tau_3})^{-\tau_1}\}^{a-1} \\ [1-\{(1+\tau_2 t^{-\tau_3})^{-\tau_1}\}^a]^{bn-1}[1-[1-\{(1+\tau_2 t^{-\tau_3})^{-\tau_1}\}^a]^b]^{m-1} \end{bmatrix}}{B(m,n)\,[I_{1-[1-\{(1+\tau_2 t^{-\tau_3})^{-\tau_1}\}^a]^b}(m,n)]}$

chrf: $H^{BKwDa}(t;a,b,m,n,\tau) = -\log[1 - I_{1-[1-\{(1+\tau_2 t^{-\tau_3})^{-\tau_1}\}^a]^b}(m,n)]$

This distribution was recently introduced by Domma and Condino (2016).

### 3.2 The $BKw$ - Singh-Maddala ($BKw-SM$) distribution

Considering the Singh-Maddala distribution (Singh and Maddala 1976) with cdf and pdf given by $G(t;\gamma) = 1 - (1+\gamma_2 t^{\gamma_3})^{-\gamma_1}$ and $g(t;\gamma) = \gamma_1\gamma_2\gamma_3 t^{(\gamma_3-1)}(1+\gamma_2 t^{\gamma_3})^{-\gamma_1-1}$ with $\gamma' = (\gamma_1, \gamma_2, \gamma_3)$ and $\gamma_i > 0$ for $i=1,2,3$ where $\gamma_2$ is a scale parameter and $\gamma_1$ and $\gamma_3$ are shape parameters. The pdf and cdf of $BKw-SM$ distribution are given by

$f^{BKwSM}(t;a,b,m,n,\gamma) = \dfrac{\begin{bmatrix} ab\gamma_1\gamma_2\gamma_3 t^{(\gamma_3-1)}(1+\gamma_2 t^{\gamma_3})^{-\gamma_1-1}\{1-(1+\gamma_2 t^{\gamma_3})^{-\gamma_1}\}^{a-1} \\ [1-\{1-(1+\gamma_2 t^{\gamma_3})^{-\gamma_1}\}^a]^{bn-1}[1-[1-\{1-(1+\gamma_2 t^{\gamma_3})^{-\gamma_1}\}^a]^b]^{m-1} \end{bmatrix}}{B(m,n)}$

and $F^{BKwSM}(t;a,b,m,n,\gamma) = I_{1-[1-\{1-(1+\gamma_2 t^{\gamma_3})^{-\gamma_1}\}^a]^b}(m,n)$

sf: $\overline{F}^{BKwSM}(t;a,b,m,n,\gamma) = 1 - I_{1-[1-\{1-(1+\gamma_2 t^{\gamma_3})^{-\gamma_1}\}^a]^b}(m,n)$

hrf: $h^{BKwSM}(t;a,b,m,n,\gamma) = \dfrac{\begin{bmatrix} ab\gamma_1\gamma_2\gamma_3 t^{(\gamma_3-1)}(1+\gamma_2 t^{\gamma_3})^{-\gamma_1-1}\{1-(1+\gamma_2 t^{\gamma_3})^{-\gamma_1}\}^{a-1} \\ [1-\{1-(1+\gamma_2 t^{\gamma_3})^{-\gamma_1}\}^a]^{bn-1}[1-[1-\{1-(1+\gamma_2 t^{\gamma_3})^{-\gamma_1}\}^a]^b]^{m-1} \end{bmatrix}}{B(m,n)\,[1 - I_{1-[1-\{1-(1+\gamma_2 t^{\gamma_3})^{-\gamma_1}\}^a]^b}(m,n)]}$

rhrf: $r^{BKwSM}(t;a,b,m,n,\gamma) = \dfrac{\begin{bmatrix} ab\gamma_1\gamma_2\gamma_3 t^{(\gamma_3-1)}(1+\gamma_2 t^{\gamma_3})^{-\gamma_1-1}\{1-(1+\gamma_2 t^{\gamma_3})^{-\gamma_1}\}^{a-1} \\ [1-\{1-(1+\gamma_2 t^{\gamma_3})^{-\gamma_1}\}^a]^{bn-1}[1-[1-\{1-(1+\gamma_2 t^{\gamma_3})^{-\gamma_1}\}^a]^b]^{m-1} \end{bmatrix}}{B(m,n)\,[I_{1-[1-\{1-(1+\gamma_2 t^{\gamma_3})^{-\gamma_1}\}^a]^b}(m,n)]}$

chrf: $H^{BKwSM}(t;a,b,m,n,\gamma) = -\log[1 - I_{1-[1-\{1-(1+\gamma_2 t^{\gamma_3})^{-\gamma_1}\}^a]^b}(m,n)]$

This distribution was recently introduced by Domma and Condino (2016).

### 3.3 The $BKw$ − exponential ($BKw-E$) distribution

Let the base line distribution be exponential with parameter $\lambda > 0$, $g(t:\lambda) = \lambda e^{-\lambda t}$, $t > 0$ and $G(t:\lambda) = 1 - e^{-\lambda t}$, $t > 0$ then for the $BKw - E$ model we get the pdf and cdf respectively as:

$$f^{BKwE}(t;m,n,a,b,\lambda) = \frac{1}{B(m,n)} a b \lambda e^{-\lambda t} \{1 - e^{-\lambda t}\}^{a-1} [1 - \{1 - e^{-\lambda t}\}^{a}]^{bn-1} [1 - [1 - \{1 - e^{-\lambda t}\}^{a}]^{b}]^{m-1}$$

and $F^{BKwE}(t;m,n,a,b,\lambda) = I_{1-[1-\{1-e^{-\lambda t}\}^{a}]^{b}}(m,n)$

sf: $\bar{F}^{BKwE}(t;m,n,a,b,\lambda) = 1 - I_{1-[1-\{1-e^{-\lambda t}\}^{a}]^{b}}(m,n)$

hrf: $h^{BKwE}(t;m,n,a,b,\lambda)$

$$= \frac{a b \lambda e^{-\lambda t} \{1 - e^{-\lambda t}\}^{a-1} [1 - \{1 - e^{-\lambda t}\}^{a}]^{bn-1} [1 - [1 - \{1 - e^{-\lambda t}\}^{a}]^{b}]^{m-1}}{B(m,n) [1 - I_{1-[1-\{1-e^{-\lambda t}\}^{a}]^{b}}(m,n)]}$$

rhrf: $r^{BKwE}(t;m,n,a,b,\lambda)$

$$= \frac{a b \lambda e^{-\lambda t} \{1 - e^{-\lambda t}\}^{a-1} [1 - \{1 - e^{-\lambda t}\}^{a}]^{bn-1} [1 - [1 - \{1 - e^{-\lambda t}\}^{a}]^{b}]^{m-1}}{B(m,n) I_{1-[1-\{1-e^{-\lambda t}\}^{a}]^{b}}(m,n)]}$$

chrf: $H^{BKwE}(t;m,n,a,b,\lambda) = -\log[1 - I_{1-[1-\{1-e^{-\lambda t}\}^{a}]^{b}}(m,n)]$

### 3.4 The $BKw - \text{Lomax}$ ($BKw - L$) distribution

Considering the Lomax distribution (Ghitany *et al.* 2007) with pdf and cdf given by $g(t:\beta,\delta) = (\beta/\delta)[1 + (t/\delta)]^{-(\beta+1)}$, $t > 0$, and $G(t:\beta,\delta) = 1 - [1 + (t/\delta)]^{-\beta}$, $\beta > 0$ and $\delta > 0$ the pdf and cdf of the $BKw - L$ distribution are given respectively by

$f^{BKwL}(t;m,n,a,b,\beta,\delta)$

$$= \frac{a b (\beta/\delta)[1 + (t/\delta)]^{-(\beta+1)} \{1 - [1 + (t/\delta)]^{-\beta}\}^{a-1} [1 - \{1 - [1 + (t/\delta)]^{-\beta}\}^{a}]^{bn-1}}{B(m,n)}$$

$$[1 - [1 - \{1 - [1 + (t/\delta)]^{-\beta}\}^{a}]^{b}]^{m-1}$$

and $F^{BKwL}(t;m,n,a,b,\beta,\delta) = I_{1-[1-\{1-[1+(t/\delta)]^{-\beta}\}^{a}]^{b}}(m,n)$

sf: $\bar{F}^{BKwL}(t;m,n,a,b,\beta,\delta) = 1 - I_{1-[1-\{1-[1+(t/\delta)]^{-\beta}\}^{a}]^{b}}(m,n)$

hrf: $h^{BKwL}(t;m,n,a,b,\beta,\delta)$

$$= \frac{a b (\beta/\delta)[1 + (t/\delta)]^{-(\beta+1)} \{1 - [1 + (t/\delta)]^{-\beta}\}^{a-1} [1 - \{1 - [1 + (t/\delta)]^{-\beta}\}^{a}]^{bn-1}}{B(m,n) [1 - I_{1-[1-\{1-[1+(t/\delta)]^{-\beta}\}^{a}]^{b}}(m,n)]}$$

$$[1 - [1 - \{1 - [1 + (t/\delta)]^{-\beta}\}^{a}]^{b}]^{m-1}$$

rhrf : $r^{BKwL}(t;m,n,a,b,\beta,\delta)$

$$= \frac{ab(\beta/\delta)[1+(t/\delta)]^{-(\beta+1)}\{1-[1+(t/\delta)]^{-\beta}\}^{a-1}[1-\{1-[1+(t/\delta)]^{-\beta}\}^a]^{bn-1}}{B(m,n)\, I_{1-[1-\{1-[1+(t/\delta)]^{-\beta}\}^a]^b}(m,n)}$$

$$[1-[1-\{1-[1+(t/\delta)]^{-\beta}\}^a]^b]^{m-1}$$

chrf: $H^{BKwL}(t;m,n,a,b,\beta,\delta) = -\log[1-I_{1-[1-\{1-[1+(t/\delta)]^{-\beta}\}^a]^b}(m,n)]$

### 3.5 The $BKw$ – Weibull ($BKw-W$) distribution

Considering the Weibull distribution (Ghitany *et al.* 2005, Zhang and Xie 2007) with parameters $\lambda>0$ and $\beta>0$ having pdf and cdf $g(t)=\lambda\beta t^{\beta-1}e^{-\lambda t^\beta}$ and $G(t)=1-e^{-\lambda t^\beta}$ respectively we get the pdf and cdf of $BKw-W$ distribution as

$$f^{BKwW}(t;m,n,a,b,\lambda,\beta)$$
$$= \frac{ab\lambda\beta t^{\beta-1}e^{-\lambda t^\beta}\{1-e^{-\lambda t^\beta}\}^{a-1}[1-\{1-e^{-\lambda t^\beta}\}^a]^{bn-1}[1-[1-\{1-e^{-\lambda t^\beta}\}^a]^b]^{m-1}}{B(m,n)}$$

and $\quad F^{BKwW}(t;m,n,a,b,\lambda,\beta) = I_{1-[1-\{1-e^{-\lambda t^\beta}\}^a]^b}(m,n)$

sf: $\quad \overline{F}^{BKwW}(t;m,n,a,b,\lambda,\beta) = 1 - I_{1-[1-\{1-e^{-\lambda t^\beta}\}^a]^b}(m,n)$

hrf: $h^{BKwW}(t;m,n,a,b,\lambda,\beta) =$

$$\frac{ab\lambda\beta t^{\beta-1}e^{-\lambda t^\beta}\{1-e^{-\lambda t^\beta}\}^{a-1}[1-\{1-e^{-\lambda t^\beta}\}^a]^{bn-1}[1-[1-\{1-e^{-\lambda t^\beta}\}^a]^b]^{m-1}}{B(m,n)\,[1-I_{1-[1-\{1-e^{-\lambda t^\beta}\}^a]^b}(m,n)]}$$

rhrf: $r^{BKwW}(t;m,n,a,b,\lambda,\beta) =$

$$\frac{ab\lambda\beta t^{\beta-1}e^{-\lambda t^\beta}\{1-e^{-\lambda t^\beta}\}^{a-1}[1-\{1-e^{-\lambda t^\beta}\}^a]^{bn-1}[1-[1-\{1-e^{-\lambda t^\beta}\}^a]^b]^{m-1}}{B(m,n)\,I_{1-[1-\{1-e^{-\lambda t^\beta}\}^a]^b}(m,n)}$$

chrf: $H^{BKwW}(t;m,n,a,b,\lambda,\beta) = -\log[1-I_{1-[1-\{1-e^{-\lambda t^\beta}\}^a]^b}(m,n)]$

### 3.6 The $BKw$ – Frechet ($BKw-Fr$) distribution

Suppose the base line distribution is the Frechet distribution (Krishna *et al.*, 2013) with pdf and cdf given by $g(t)=\lambda\delta^\lambda t^{-(\lambda+1)}e^{-(\delta/t)^\lambda}$ and $G(t)=e^{-(\delta/t)^\lambda}$, $t>0$ respectively, then the corresponding pdf and cdf of $BKw-Fr$ distribution becomes

$$f^{BKwFr}(t;m,n,a,b,\lambda,\delta)$$
$$= \frac{ab\lambda\delta^\lambda t^{-(\lambda+1)}e^{-(\delta/t)^\lambda}\{e^{-(\delta/t)^\lambda}\}^{a-1}[1-\{e^{-(\delta/t)^\lambda}\}^a]^{bn-1}[1-[1-\{e^{-(\delta/t)^\lambda}\}^a]^b]^{m-1}}{B(m,n)}$$

and $F^{BKwFr}(t;m,n,a,b,\lambda,\delta) = I_{1-[1-\{e^{-(\delta/t)^\lambda}\}^a]^b}(m,n)$

sf: $\quad \overline{F}^{BKwFr}(t;m,n,a,b,\lambda,\delta) = 1 - I_{1-[1-\{e^{-(\delta/t)^\lambda}\}^a]^b}(m,n)$

hrf: $h^{BKwFr}(t;m,n,a,b,\lambda,\delta) =$

$$\frac{ab\lambda\delta^{\lambda} t^{-(\lambda+1)} e^{-(\delta/t)^{\lambda}} \{e^{-(\delta/t)^{\lambda}}\}^{a-1} [1-\{e^{-(\delta/t)^{\lambda}}\}^{a}]^{bn-1} [1-[1-\{e^{-(\delta/t)^{\lambda}}\}^{a}]^{b}]^{m-1}}{B(m,n) \ [1-I_{1-[1-\{e^{-(\delta/t)^{\lambda}}\}^{a}]^{b}}(m,n)]}$$

rhrf: $r^{BKwFr}(t;m,n,a,b,\lambda,\delta) =$

$$\frac{ab\lambda\delta^{\lambda} t^{-(\lambda+1)} e^{-(\delta/t)^{\lambda}} \{e^{-(\delta/t)^{\lambda}}\}^{a-1} [1-\{e^{-(\delta/t)^{\lambda}}\}^{a}]^{bn-1} [1-[1-\{e^{-(\delta/t)^{\lambda}}\}^{a}]^{b}]^{m-1}}{B(m,n) \ [I_{1-[1-\{e^{-(\delta/t)^{\lambda}}\}^{a}]^{b}}(m,n)]}$$

chrf: $H^{BKwFr}(t;m,n,a,b,\lambda,\delta) = -\log [1-I_{1-[1-\{e^{-(\delta/t)^{\lambda}}\}^{a}]^{b}}(m,n)]$

### 3.7 The $BKw$ – Gompertz ($BKw-Go$) distribution

Next by taking the Gompertz distribution (Gieser et al. 1998) with pdf and cdf $g(t) = \beta e^{\lambda t} e^{-\frac{\beta}{\lambda}(e^{\lambda t}-1)}$ and $G(t) = 1 - e^{-\frac{\beta}{\lambda}(e^{\lambda t}-1)}$, $\beta > 0, \lambda > 0, t > 0$ respectively, we get the pdf and cdf of $BKw-Go$ distribution as

$f^{BKwGo}(t;m,n,a,b,\lambda,\beta)$

$$= \frac{ab\beta e^{\lambda t} e^{-\frac{\beta}{\lambda}(e^{\lambda t}-1)} \{1-e^{-\frac{\beta}{\lambda}(e^{\lambda t}-1)}\}^{a-1} [1-\{1-e^{-\frac{\beta}{\lambda}(e^{\lambda t}-1)}\}^{a}]^{bn-1} [1-[1-\{1-e^{-\frac{\beta}{\lambda}(e^{\lambda t}-1)}\}^{a}]^{b}]^{m-1}}{B(m,n)} \text{ and}$$

$F^{BKwGo}(t;m,n,a,b,\lambda,\beta) = I_{1-[1-\{1-e^{-\frac{\beta}{\lambda}(e^{\lambda t}-1)}\}^{a}]^{b}}(m,n)$

sf: $\overline{F}^{BKwG}(t;m,n,a,b,\lambda,\beta) = 1 - I_{1-[1-\{1-e^{-\frac{\beta}{\lambda}(e^{\lambda t}-1)}\}^{a}]^{b}}(m,n)$

hrf: $h^{BKwG}(t;m,n,a,b,\lambda,\beta) =$

$$\frac{ab\beta e^{\lambda t} e^{-\frac{\beta}{\lambda}(e^{\lambda t}-1)} \{1-e^{-\frac{\beta}{\lambda}(e^{\lambda t}-1)}\}^{a-1} [1-\{1-e^{-\frac{\beta}{\lambda}(e^{\lambda t}-1)}\}^{a}]^{bn-1} [1-[1-\{1-e^{-\frac{\beta}{\lambda}(e^{\lambda t}-1)}\}^{a}]^{b}]^{m-1}}{B(m,n) \ [1-I_{1-[1-\{1-e^{-\frac{\beta}{\lambda}(e^{\lambda t}-1)}\}^{a}]^{b}}(m,n)]}$$

rhrf: $r^{BKwG}(t;m,n,a,b,\lambda,\beta) =$

$$\frac{ab\beta e^{\lambda t} e^{-\frac{\beta}{\lambda}(e^{\lambda t}-1)} \{1-e^{-\frac{\beta}{\lambda}(e^{\lambda t}-1)}\}^{a-1} [1-\{1-e^{-\frac{\beta}{\lambda}(e^{\lambda t}-1)}\}^{a}]^{bn-1} [1-[1-\{1-e^{-\frac{\beta}{\lambda}(e^{\lambda t}-1)}\}^{a}]^{b}]^{m-1}}{B(m,n) \ I_{1-[1-\{1-e^{-\frac{\beta}{\lambda}(e^{\lambda t}-1)}\}^{a}]^{b}}(m,n)}$$

chrf: $H^{BKwG}(t;m,n,a,b,\lambda,\beta) = -\log [1-I_{1-[1-\{1-e^{-\frac{\beta}{\lambda}(e^{\lambda t}-1)}\}^{a}]^{b}}(m,n)]$

### 3.8 The $BKw$ – Extended Weibull ($BKw-EW$) distribution

The pdf and the cdf of the extended Weibull (*EW*) distributions of (Gurvich *et al.* 1997) is given by

$g(t:\delta,\vartheta) = \delta \exp[-\delta Z(t:\vartheta)]z(t:\vartheta)$ and $G(t:\delta,\vartheta) = 1 - \exp[-\delta Z(t:\vartheta)]$, $\quad t \in D \subseteq R_{+}, \delta > 0$

where $Z(t:\vartheta)$ is a non-negative monotonically increasing function which depends on the parameter vector $\vartheta$ and $z(t:\vartheta)$ is the derivative of $Z(t:\vartheta)$.

By considering *EW* as the base line distribution we derive pdf and cdf of the $BKw-EW$ as

$$f^{BKwEW}(t;m,n,a,b,\delta,\vartheta)$$
$$= \frac{ab\delta \exp[-\delta Z(t:\vartheta)]z(t:\vartheta)\{1-\exp[-\delta Z(t:\vartheta)]\}^{a-1}[1-\{1-\exp[-\delta Z(t:\vartheta)]\}^{a}]^{bn-1}}{B(m,n)}$$
$$\times [1-[1-\{1-\exp[-\delta Z(t:\vartheta)]\}^{a}]^{b}]^{m-1}$$

and $F^{BKwEW}(t;m,n,a,b,\delta,\vartheta) = I_{1-[1-\{1-\exp[-\delta Z(t:\vartheta)]\}^{a}]^{b}}(m,n)$

Important models can be seen as particular cases with different choices of $Z(t:\vartheta)$:

(i) $Z(t:\vartheta) = t$: exponential distribution.

(ii) $Z(t:\vartheta) = t^2$: Rayleigh (Burr type-X) distribution.

(iii) $Z(t:\vartheta) = \log(t/k)$: Pareto distribution

(iv) $Z(t:\vartheta) = \beta^{-1}[\exp(\beta t) - 1]$: Gompertz distribution.

sf: $\overline{F}^{BKwEW}(t;m,n,a,b,\delta,\vartheta) = 1 - I_{1-[1-\{1-\exp[-\delta Z(t:\vartheta)]\}^{a}]^{b}}(m,n)$

hrf: $h^{BKwEW}(t;m,n,a,b,\delta,\vartheta) =$

$$\frac{ab\delta \exp[-\delta Z(t:\vartheta)]z(t:\vartheta)\{1-\exp[-\delta Z(t:\vartheta)]\}^{a-1}[1-\{1-\exp[-\delta Z(t:\vartheta)]\}^{a}]^{bn-1}}{B(m,n)\ [1-I_{1-[1-\{1-\exp[-\delta Z(t:\vartheta)]\}^{a}]^{b}}(m,n)]}$$
$$\times [1-[1-\{1-\exp[-\delta Z(t:\vartheta)]\}^{a}]^{b}]^{m-1}$$

rhrf: $r^{BKwEW}(t;m,n,a,b,\delta,\vartheta) =$

$$\frac{ab\delta \exp[-\delta Z(t:\vartheta)]z(t:\vartheta)\{1-\exp[-\delta Z(t:\vartheta)]\}^{a-1}[1-\{1-\exp[-\delta Z(t:\vartheta)]\}^{a}]^{bn-1}}{B(m,n)\ I_{1-[1-\{1-\exp[-\delta Z(t:\vartheta)]\}^{a}]^{b}}(m,n)}$$
$$\times [1-[1-\{1-\exp[-\delta Z(t:\vartheta)]\}^{a}]^{b}]^{m-1}$$

chrf: $H^{BKwEW}(t;m,n,a,b,\delta,\vartheta) = -\log[1-I_{1-[1-\{1-\exp[-\delta Z(t:\vartheta)]\}^{a}]^{b}}(m,n)]$

**3.9 The $BKw$–Extended Modified Weibull ($BKw-EMW$) distribution**

Considering the modified Weibull (MW) distribution (Sarhan and Zaindin 2013) given by cdf and pdf

$$G(t;\sigma,\beta,\gamma) = 1 - \exp[-\sigma t - \beta t^{\gamma}], \ t > 0, \gamma > 0, \sigma, \beta \geq 0, \sigma + \beta > 0$$

and $g(t;\sigma,\beta,\gamma) = (\sigma + \beta\gamma t^{\gamma-1})\exp[-\sigma t - \beta t^{\gamma}]$ respectively the corresponding pdf and cdf of $BKw-EMW$ are given by

$f^{BKwEMW}(t;m,n,a,b,\sigma,\beta,\gamma)$

$$= \frac{\begin{bmatrix} ab(\sigma + \beta\gamma t^{\gamma-1})\exp[-\sigma t - \beta t^\gamma]\{1-\exp[-\sigma t - \beta t^\gamma]\}^{a-1} \\ [1-\{1-\exp[-\sigma t - \beta t^\gamma]\}^a]^{bn-1}[1-[1-\{1-\exp[-\sigma t - \beta t^\gamma]\}^a]^b]^{m-1} \end{bmatrix}}{B(m,n)}$$

and $F^{BKwEMW}(t;m,n,a,b,\sigma,\beta,\gamma) = I_{1-[1-\{1-\exp[-\sigma t-\beta t^\gamma]\}^a]^b}(m,n)$.

sf: $\overline{F}^{BKwEMW}(t;m,n,a,b,\sigma,\beta,\gamma) = 1 - I_{1-[1-\{1-\exp[-\sigma t-\beta t^\gamma]\}^a]^b}(m,n)$

hrf: $h^{BKwEMW}(t;m,n,a,b,\sigma,\beta,\gamma) =$

$$\frac{ab(\sigma+\beta\gamma t^{\gamma-1})\exp[-\sigma t - \beta t^\gamma]\{1-\exp[-\sigma t - \beta t^\gamma]\}^{a-1}[1-\{1-\exp[-\sigma t - \beta t^\gamma]\}^a]^{bn-1}}{B(m,n)\ [1-I_{1-[1-\{1-\exp[-\sigma t-\beta t^\gamma]\}^a]^b}(m,n)]}$$

$$\times [1-[1-\{1-\exp[-\sigma t - \beta t^\gamma]\}^a]^b]^{m-1}$$

rhrf: $r^{BKwEMW}(t;m,n,a,b,\sigma,\beta,\gamma) =$

$$\frac{ab(\sigma+\beta\gamma t^{\gamma-1})\exp[-\sigma t - \beta t^\gamma]\{1-\exp[-\sigma t - \beta t^\gamma]\}^{a-1}[1-\{1-\exp[-\sigma t - \beta t^\gamma]\}^a]^{bn-1}}{B(m,n)\ I_{1-[1-\{1-\exp[-\sigma t-\beta t^\gamma]\}^a]^b}(m,n)}$$

$$\times [1-[1-\{1-\exp[-\sigma t - \beta t^\gamma]\}^a]^b]^{m-1}$$

chrf: $H^{BKwEMW}(t;m,n,a,b,\sigma,\beta,\gamma) = -\log[1 - I_{1-[1-\{1-\exp[-\sigma t-\beta t^\gamma]\}^a]^b}(m,n)]$

### 3.10 The $BKw$ – Extended Exponentiated Pareto ($BKw$ – $EEP$) distribution

The pdf and cdf of the exponentiated Pareto distribution, of Nadarajah (2005), are given respectively by $g(t) = \gamma k \theta^k t^{-(k+1)}[1-(\theta/t)^k]^{\gamma-1}$ and $G(t) = [1-(\theta/t)^k]^\gamma$, $x > \theta$ and $\theta, k, \gamma > 0$.

Thus the pdf and the cdf of $BKw - EEP$ distribution are given by

$f^{BKwEEP}(t;m,n,a,b,\theta,k,\gamma)$

$$= \frac{\begin{bmatrix} ab\gamma k \theta^k t^{-(k+1)}[1-(\theta/t)^k]^{\gamma-1}\{[1-(\theta/t)^k]^\gamma\}^{a-1}[1-\{[1-(\theta/t)^k]^\gamma\}^a]^{bn-1} \\ [1-[1-\{[1-(\theta/t)^k]^\gamma\}^a]^b]^{m-1} \end{bmatrix}}{B(m,n)}$$

and $F^{BKwEEP}(t;m,n,a,b,\theta,k,\gamma) = I_{1-[1-\{[1-(\theta/t)^k]^\gamma\}^a]^b}(m,n)$.

sf: $\overline{F}^{BKwEEP}(t;m,n,a,b,\theta,k,\gamma) = 1 - I_{1-[1-\{[1-(\theta/t)^k]^\gamma\}^a]^b}(m,n)$

hrf: $h^{BKwEEP}(t;m,n,a,b,\theta,k,\gamma) =$

$$\frac{ab\gamma k \theta^k t^{-(k+1)}[1-(\theta/t)^k]^{\gamma-1}\{[1-(\theta/t)^k]^\gamma\}^{a-1}[1-\{[1-(\theta/t)^k]^\gamma\}^a]^{bn-1}[1-[1-\{[1-(\theta/t)^k]^\gamma\}^a]^b]^{m-1}}{B(m,n)\ 1-I_{1-[1-\{[1-(\theta/t)^k]^\gamma\}^a]^b}(m,n)}$$

rhrf: $r^{BKwEEP}(t;m,n,a,b,\theta,k,\gamma) =$

$$\frac{ab\gamma k \theta^k t^{-(k+1)}[1-(\theta/t)^k]^{\gamma-1}\{[1-(\theta/t)^k]^\gamma\}^{a-1}[1-\{[1-(\theta/t)^k]^\gamma\}^a]^{bn-1}[1-[1-\{[1-(\theta/t)^k]^\gamma\}^a]^b]^{m-1}}{B(m,n)\ I_{1-[1-\{[1-(\theta/t)^k]^\gamma\}^a]^b}(m,n)}$$

chrf: $H^{BKwEEP}(t;m,n,a,b,\theta,k,\gamma) = -\log[1 - I_{1-[1-\{[1-(\theta/t)^k]^\gamma\}^a]^b}(m,n)]$

## 4. General results for the Beta Kumaraswamy-G ($BKw-G$) family of distributions

In this section we derive some general results for the proposed $BKw-G$ family.

### 4.1 Expansions of pdf and cdf

By using binomial expansion in (6), we obtain

$$f^{BKwG}(t;a,b,m,n) = \frac{1}{B(m,n)} a b g(t) G(t)^{a-1} [1-G(t)^a]^{bn-1} [1-[1-G(t)^a]^b]^{(m-1)}$$

$$= \frac{1}{B(m,n)} a b g(t) G(t)^{a-1} [1-G(t)^a]^{bn-1} \sum_{j=0}^{m-1} \binom{m-1}{j} (-1)^j [1-G(t)^a]^{bj} \quad (10)$$

$$= \sum_{j=0}^{m-1} \beta'_j f^{KwG}(t;a,b(j+n)) \quad (11)$$

Where $\beta'_j = \frac{(-1)^{j+1}}{B(m,n)(j+n)} \binom{m-1}{j}$. Adjusting further we get from (10)

$$= \sum_{j=0}^{m-1} \beta'_j \frac{d}{dt} [\overline{F}^{KwG}(t;a,b(j+n))]$$

$$= \sum_{j=0}^{m-1} \beta'_j \frac{d}{dt} [\overline{F}^{KwG}(t;a,b)]^{j+n} \quad (12)$$

$$= f^{KwG}(t;a,b) \sum_{j=0}^{m-1} \{\beta'_j (j+n)\} [\overline{F}^{KwG}(t;a,b)]^{j+n-1}$$

$$= f^{KwG}(t;a,b) \sum_{j=0}^{m-1} \beta_j [\overline{F}^{KwG}(t;a,b)]^{j+n-1} \quad (13)$$

Where $\beta_j = \beta'_j (j+n)$

Alternatively, we can expand the pdf as

$$f^{BKwG}(t;a,b,m,n) = f^{KwG}(t;a,b) \sum_{j=0}^{m-1} \beta_j [\overline{F}^{KwG}(t;a,b)]^{j+n-1}$$

$$= f^{KwG}(t;a,b) \sum_{j=0}^{m-1} \beta_j [1-F^{KwG}(t;a,b)]^{j+n-1}$$

$$= f^{KwG}(t;a,b) \sum_{j=0}^{m-1} \beta_j \sum_{l=0}^{j+n-1} \binom{j+n-1}{l} (-1)^l [F^{KwG}(t;a,b)]^l$$

$$= f^{KwG}(t;a,b) \sum_{l=0}^{j+n-1} \eta_l [F^{KwG}(t;a,b)]^l \quad (14)$$

Where $\eta_l = \sum_{j=0}^{m-1} (-1)^l \beta_j \binom{j+n-1}{l}$.

We can expand the cdf with the following result

$$I_z(a,b) = \frac{z^a}{B(m,n)} \sum_{n=0}^{\infty} \binom{b-1}{n} \frac{(-1)^n z^n}{(a+n)} \quad (15)$$

Where $(x)_n$ is a Pochhammer symbol (See "Incomplete Beta Function" *From Math World*-A Wolfram Web Resource. http://mathworld. Wolfram.com/Incomplete Beta Function. html).

Using (15) in (7) we have

$$F^{BKwG}(t;a,b,m,n) = \frac{1}{B(m,n)} \left(1-[1-G(t)^a]^b\right)^m \sum_{i=0}^{\infty} \binom{n-1}{i} \frac{(-1)^i}{(m+i)} \left(1-[1-G(t)^a]^b\right)^i$$

$$= \sum_{i=0}^{\infty} \frac{(-1)^i}{B(m,n)(m+i)} \binom{n-1}{i} [F^{KwG}(t;a,b)]^{m+i}$$

$$= \sum_{i=0}^{\infty} \frac{(-1)^i}{B(m,n)(m+i)} \binom{n-1}{i} \sum_{j=0}^{\infty} \binom{m+i}{j} (-1)^j [\overline{F}^{KwG}(t;a,b)]^j$$

$$= \sum_{i,j=0}^{\infty} \frac{(-1)^{i+j}}{B(m,n)(m+i)} \binom{n-1}{i}\binom{m+i}{j} \sum_{k=0}^{j} \binom{j}{k} (-1)^k [F^{KwG}(t;a,b)]^k$$

$$= \sum_{i,j=0}^{\infty} \sum_{k=0}^{j} \frac{(-1)^{i+j+k}}{B(m,n)(m+i)} \binom{n-1}{i}\binom{m+i}{j}\binom{j}{k} [F^{KwG}(t;a,b)]^k$$

By exchanging the indices $j$ and $k$ in the sum symbol, we get

$$F^{BKwG}(t;a,b,m,n) = \sum_{i,k=0}^{\infty} \sum_{j=k}^{\infty} \frac{(-1)^{i+j+k}}{B(m,n)(m+i)} \binom{n-1}{i}\binom{m+i}{j}\binom{j}{k} [F^{KwG}(t;a,b)]^k$$

and then $F^{BKwG}(t;a,b,m,n) = \sum_{k=0}^{\infty} \mu_k F^{KwG}(t;a,b)^k \quad (16)$

Where $\mu_k = \sum_{i=0}^{\infty} \sum_{j=k}^{\infty} \frac{(-1)^{i+j+k}}{B(m,n)(m+i)} \binom{n-1}{i}\binom{m+i}{j}\binom{j}{k}$

Similarly, an expansion for the cdf of $BKw-G$ can be derives as (David Osborn & Richard Madey 1968)

$$F^{BKwG}(t;a,b,m,n) = I_{1-[1-G(t)^a]^b}(m,n) = I_{F^{KwG}(t;a,b)}(m,n)$$

$$= \sum_{p=m}^{m+n-1} \binom{m+n-1}{p} [F^{KwG}(t;a,b)]^p [\overline{F}^{KwG}(t;a,b)]^{m+n-1-p}$$

$$= \sum_{p=m}^{m+n-1} \binom{m+n-1}{p} [F^{KwG}(t;a,b)]^p \sum_{q=0}^{m+n-1-p} \binom{m+n-1-p}{q} [-F^{KwG}(t;a,b)]^q$$

$$= \sum_{p=m}^{m+n-1} \sum_{q=0}^{m+n-1-p} \binom{m+n-1}{p}\binom{m+n-1-p}{q} (-1)^q [F^{KwG}(t;a,b)]^{p+q}$$

$$= \sum_{p=m}^{m+n-1} \sum_{q=0}^{m+n-1-p} \lambda_{p,q} \, [F^{KwG}(t;a,b)]^{p+q}.$$

$$\text{Where } \lambda_{p,q} = \binom{m+n-1}{p} \binom{m+n-1-p}{q} (-1)^q.$$

## 4.2 Order statistics

Suppose $T_1, T_2, ..., T_n$ is a random sample from any $BKw-G$ distribution. Let $T_{r:n}$ denote the $r^{th}$ order statistics. The pdf of $T_{r:n}$ can be expressed as

$$f_{r:n}(t;a,b,m,n) = \frac{n!}{(r-1)!(n-r)!} f^{BKwG}(t) \, F^{BKwG}(t)^{r-1} \{1 - F^{BKwG}(t)\}^{n-r}$$

$$= \frac{n!}{(r-1)!(n-r)!} \sum_{j=0}^{n-r} (-1)^j \binom{n-r}{j} f^{BKwG}(t) \, F^{BKwG}(t)^{j+r-1}$$

Now using the general expansion of the pdf and cdf of the $BKw-G$ distribution we get the pdf of the $r^{th}$ order statistics for of the $BKw-G$ as

$$f_{r:n}(t;a,b,m,n) = \frac{n!}{(r-1)!(n-r)!} \sum_{j=0}^{n-r} (-1)^j \binom{n-r}{j} \left\{ f^{KwG}(t;a,b) \sum_{l=0}^{j+n-1} \eta_l [F^{KwG}(t;a,b)]^l \right\}$$

$$\left\{ \sum_{k=0}^{\infty} \mu_k [F^{KwG}(t;a,b)]^k \right\}^{j+r-1}$$

Where $\eta_l$ and $\mu_k$ defined in section 4.

Now $\left\{ \sum_{k=0}^{\infty} \mu_k [F^{KwG}(t;a,b)]^k \right\}^{j+r-1} = \sum_{k=0}^{\infty} d_{j+r-1,k} [F^{KwG}(t;a,b)]^k$

Where $d_{j+r-1,k} = \frac{1}{k\mu_0} \sum_{c=1}^{k} [c(j+r)-k] \mu_c \, d_{j+r-1,\, k-c}$ (Nadarajah et. al 2015)

Therefore the density function of the $r^{th}$ order statistics of $BKw-G$ distribution can be expressed as

$$f_{r:n}(t;a,b,m,n) = \frac{n!}{(r-1)!(n-r)!} \sum_{j=0}^{n-r} (-1)^j \binom{n-r}{j} \left\{ f^{KwG}(t;a,b) \sum_{l=0}^{j+n-1} \eta_l [F^{KwG}(t;a,b)]^l \right\}$$

$$\left( \sum_{k=0}^{\infty} d_{j+r-1,k} [F^{KwG}(t;a,b)]^k \right)$$

$$= \frac{n!}{(r-1)!(n-r)!} \sum_{j=0}^{n-r} (-1)^j \binom{n-r}{j} \left\{ f^{KwG}(t;a,b) \sum_{l=0}^{j+n-1} \sum_{k=0}^{\infty} \eta_l \, d_{j+r-1,k} [F^{KwG}(t;a,b)]^{k+l} \right\}$$

$$= f^{KwG}(t;a,b) \sum_{l=0}^{j+n-1} \sum_{k=0}^{\infty} \gamma_{l,k} [F^{KwG}(t;a,b)]^{k+l}$$

$$= f^{KwG}(t;a,b) \sum_{l=0}^{j+n-1} \sum_{k=0}^{\infty} \gamma_{l,k} \sum_{z=0}^{k+l} \binom{k+l}{z} (-1)^z [\overline{F}^{KwG}(t;a,b)]^z$$

$$= f^{KwG}(t;a,b) \sum_{z=0}^{k+l} \chi_z [\overline{F}^{KwG}(t;a,b)]^z \qquad (17)$$

$$= \sum_{z=0}^{k+l} \frac{\chi_z}{z+1} \frac{d}{dt} [\overline{F}^{KwG}(t;a,b)]^{z+1} = \sum_{z=0}^{k+l} \chi'_z f^{KwG}(t;a,b(z+1))$$

Where $\chi'_z = \sum_{l=0}^{j+n-1} \sum_{k=0}^{\infty} \frac{(-1)^{z+1} \gamma_{l,k}}{z+1} \binom{k+l}{z}$, $\chi_z = \chi'_z(z+1)$

$\gamma_{l,k} = \frac{n!}{(r-1)!(n-r)!} \sum_{j=0}^{n-r} (-1)^j \binom{n-r}{j} \eta_l \, d_{j+r-1,k}$ and $\eta_l$, $d_{j+r-1,k}$ are defined above.

**4.3 Probability weighted moments**

The probability weighted moments (PWMs), first proposed by Greenwood et al. (1979), are expectations of certain functions of a random variable whose mean exists. The $(p,q,r)^{th}$ PWM of $T$ is defined by $\Gamma_{p,q,r} = \int_{-\infty}^{\infty} t^p F(t)^q [1-F(t)]^r f(t) dt$

From equations (13) and (14) the $s^{th}$ moment of $T$ can be written either as

$$E(T^s) = \int_{-\infty}^{\infty} t^s \, f^{BKwG}(t;a,b,m,n) \, dt$$

$$= \sum_{j=0}^{m-1} \beta_j \int_{-\infty}^{\infty} t^s \, [\overline{F}^{KwG}(t;a,b)]^{j+n-1} f^{KwG}(t;a,b) \, dt$$

$$= \sum_{j=0}^{m-1} \beta_j \int_{-\infty}^{\infty} t^s \, [\{1-G(t)^a\}^b]^{j+n-1} \, ab\, g(t) G(t)^{a-1} [1-G(t)^a]^{b-1} \, dt$$

$$= \sum_{j=0}^{m-1} \beta_j \, \Gamma_{s,0,j+n-1}$$

or $E(T^s) = \sum_{l=0}^{j+n-1} \eta_l \int_{-\infty}^{\infty} t^s \, [F^{KwG}(t;a,b)]^l f^{KwG}(t;a,b) \, dt$

$$= \sum_{l=0}^{j+n-1} \eta_l \int_{-\infty}^{\infty} t^s \, [1-\{1-G(t)^a\}^b]^l \, ab\, g(t) G(t)^{a-1} [1-G(t)^a]^{b-1} \, dt$$

$$= \sum_{l=0}^{j+n-1} \eta_l \, \Gamma_{s,l,0}$$

Where $\Gamma_{p,q,r} = \int_{-\infty}^{\infty} t^p \, [1-\{1-G(t)^a\}^b]^q [\{1-G(t)^a\}^b]^r \, ab\, g(t) G(t)^{a-1} [1-G(t)^a]^{b-1} \, dt$

is the PWM of $Kw\text{-}G(a,b)$ distribution.

Therefore the moments of the $BKw\text{-}G(m,n,a,b)$ can be expresses in terms of the PWMs of $Kw\text{-}G(a,b)$ (Cordeiro and de Castro, 2011). The PWM method can generally be used for estimating parameters quantiles of generalized distributions. These moments have low variance and no severe biases, and they compare favourably with estimators obtained by maximum likelihood.

Proceeding as above we can derive $s^{th}$ moment of the $r^{th}$ order statistic $T_{r:n}$, in a random sample of size $n$ from $BKw-G$ on using equation (17) as

$$E(T^s_{r;n}) = \sum_{z=0}^{k+l} \chi_z \Gamma_{s,0,z}$$

Where $\beta_j, \eta_l$ and $\chi_z$ defined in above.

## 4.4 Moment generating function

The moment generating function of $BKw-G$ family can be easily expressed in terms of those of the exponentiated $Kw-G$ (Cordeiro and de Castro, 2011) distribution using the results of section 4. For example using equation (12) it can be seen that

$$M_T(s) = E[e^{sT}] = \int_{-\infty}^{\infty} e^{st} f(t;m,n,a,b)\,dt = \int_{-\infty}^{\infty} e^{st} \sum_{j=0}^{m-1} \beta'_j \frac{d}{dt}[\overline{F}^{KwG}(t;a,b)]^{j+n}\,dt$$

$$= \sum_{j=0}^{m-1} \beta'_j \int_{-\infty}^{\infty} e^{st} \frac{d}{dt}\{\overline{F}^{KwG}(t;a,b)\}^{j+n}\,dt = \sum_{j=0}^{m-1} \beta'_j M_X(s)$$

Where $M_X(s)$ is the mgf of a $Kw-G$ (Cordeiro and de Castro, 2011) distribution.

## 4.5 Renyi entropy

The entropy of a random variable is a measure of uncertainty variation and has been used in various situations in science and engineering. The Rényi entropy is defined by

$$I_R(\delta) = (1-\delta)^{-1} \log\left(\int_{-\infty}^{\infty} f(t)^\delta\,dt\right)$$

where $\delta > 0$ and $\delta \neq 1$ For furthers details, see Song (2001). Using binomial expansion in (6) we can write

$$f^{BKwG}(t;a,b,m,n)^\delta = [\frac{1}{B(m,n)} a b\, g(t) G(t)^{a-1}[1-G(t)^a]^{bn-1}]^\delta\,[1-[1-G(t)^a]^b]^{(m-1)\delta}$$

$$= \frac{1}{B(m,n)^\delta}[a b\, g(t) G(t)^{a-1}[1-G(t)^a]^{bn-1}]^\delta \sum_{\alpha=0}^{(m-1)\delta} \binom{(m-1)\delta}{\alpha}(-1)^\alpha [1-G(t)^a]^{b\alpha}$$

Thus the Rényi entropy of $T$ can be obtained as

$$I_R(\delta) = (1-\delta)^{-1} \log\left( \sum_{\alpha=0}^{(m-1)\delta} R_\alpha \int_{-\infty}^{\infty} [a b g(t) G(t)^{a-1} [1-G(t)^a]^{bn-1}]^\delta [1-G(t)^a]^{b\alpha} dt \right)$$

$$= (1-\delta)^{-1} \log\left( \sum_{\alpha=0}^{(m-1)\delta} R_\alpha \int_{-\infty}^{\infty} [f^{KwG}(t;a,bn)]^\delta [\overline{F}^{KwG}(t;a,b)]^\alpha dt \right)$$

Where $R_\alpha = \dfrac{1}{B(m,n)^\delta} \dbinom{(m-1)\delta}{\alpha} (-1)^\alpha$.

## 4.6 Quantile power series and random sample generation

The quantile function of $T$, $t = Q(u) = F^{-1}(u)$, can be obtained by inverting (7). Let $z = Q_{m,n}(u)$ be the beta quantile function. Then,

$$F^{BKwG}(t;a,b,m,n) = I_{1-[1-G(t)^a]^b}(m,n) = u$$

$$\Rightarrow 1 - \{1 - G(t)^a\}^b = Q_{m,n}(u)$$

$$\Rightarrow \{1 - G(t)^a\}^b = 1 - Q_{m,n}(u)$$

$$\Rightarrow G(t)^a = 1 - \{1 - Q_{m,n}(u)\}^{1/b}$$

$$\Rightarrow G(t) = [1 - \{1 - Q_{m,n}(u)\}^{1/b}]^{1/a}$$

$$\Rightarrow t = Q_G[[1 - \{1 - Q_{m,n}(u)\}^{1/b}]^{1/a}]$$

It is possible to obtain an expansion for $Q_{m,n}(u)$ as

$$z = Q_{m,n}(u) = \sum_{i=0}^{\infty} e_i u^{i/m}$$

(see "Power series" From MathWorld-A Wolfram Web Resource. http:// mathworld. wolfram.com/ PowerSeries.html)

Where $e_i = [m B(m,n)]^{1/m} d_i$ and $d_0 = 0$, $d_1 = 1$, $d_2 = (n-1)/(m+1)$,

$$d_3 = \frac{(n-1)(m^2 + 3mn - m + 5n - 4)}{2(m+1)^2(m+2)}$$

$$d_4 = (n-1)[m^4 + (6n-1)m^3 + (n+2)(8n-5)m^2 + (33n^2 - 30n + 4)m + n(31n-47) + 18] / [3(m+1)^3(m+2)(m+3)]\ldots$$

The Bowley skewness (Kenney and Keeping 1962) measures and Moors kurtosis (Moors 1988) measure are robust and less sensitive to outliers and exist even for distributions without moments. For $BKw - G$ family these measures are given by

$$B = \frac{Q(3/4) + Q(1/4) - 2Q(1/2)}{Q(3/4) - Q(1/4)} \quad \text{and} \quad M = \frac{Q(3/8) - Q(1/8) + Q(7/8) - Q(5/8)}{Q(6/8) - Q(2/8)}$$

For example, let the $G$ be exponential distribution with parameter $\lambda > 0$, having pdf and cdf as $g(t:\lambda) = \lambda e^{-\lambda t}$, $t > 0$ and $G(t:\lambda) = 1 - e^{-\lambda t}$, respectively. Then the $p^{th}$ quantile is obtained as $-(1/\lambda)\log[1-p]$. Therefore, the $p^{th}$ quantile $t_p$, of $BKw-E$ is given by

$$t_p = -\frac{1}{\lambda}\log\left[1 - [\{1 - \{1 - Q_{m,n}(p)\}^{1/b}\}^{1/a}]\right]$$

**4.7 Asymptotes and shapes**

Here we investigate the asymptotic shapes of the proposed family following the method followed in Alizadeh *et al.*, (2015).

**Proposition 1.** The asymptotes of equations (6), (7) and (8) as $t \to 0$ are given by

$$f(t;a,b,m,n) \sim \frac{a\,b\,g(t)G(t)^{a-1}[1-[1-G(t)^a]^b]^{m-1}}{B(m,n)} \quad \text{as } G(t) \to 0$$

$$F(t;a,b,m,n) \sim \frac{\{1-[1-G(t)^a]^b\}^m}{B(m,n)\,m} \quad \text{as } G(t) \to 0$$

$$h(t;a,b,m,n) \sim \frac{a\,b\,g(t)G(t)^{a-1}[1-[1-G(t)^a]^b]^{m-1}}{B(m,n)} \quad \text{as } G(t) \to 0$$

**Proposition 2.** The asymptotes of equations (6), (7) and (8) as $t \to \infty$ are given by

$$f(t;a,b,m,n) \sim \frac{a\,b\,g(t)[1-G(t)^a]^{bn-1}}{B(m,n)} \quad \text{as } t \to \infty$$

$$1 - F(t;a,b,m,n) \sim \frac{[1-G(t)^a]^{bn}}{n\,B(m,n)} \quad \text{as } t \to \infty$$

$$h(t;a,b,m,n) \sim a\,b\,n\,g(t)[1-G(t)^a]^{-1} \quad \text{as } t \to \infty$$

The shapes of the density and hazard rate functions can be described analytically. The critical points of the $BKw - G$ density function are the roots of the equation:

$$\frac{d\log[f(t;a,b,m,n)]}{dt}$$

$$= \frac{g'(t)}{g(t)} + (a-1)\frac{g(t)}{G(t)} + a(1-bn)\frac{g(t)G(t)^{a-1}}{1-G(t)^a} + (m-1)\frac{a\,b\,g(t)G(t)^{a-1}[1-G(t)^a]^{b-1}}{1-[1-G(t)^a]^b} \quad (18)$$

There may be more than one root to (18). If $t = t_0$ is a root of (18) then it corresponds to a local maximum, a local minimum or a point of inflexion depending on whether $\kappa(t_0) < 0$, $\kappa(t_0) > 0$ *or* $\kappa(t_0) = 0$ where

$$\kappa(t) = \frac{d^2}{dt^2}\log[f(t;a,b,m,n)]$$

$$= \frac{g(t)g''(t) - [g'(t)]^2}{g(t)^2} + (a-1)\frac{G(t)g'(t) - g(t)^2}{G(t)^2}$$

$$+ a(1-bn)\left[\frac{g'(t)G(t)^{a-1}}{1-G(t)^a} + \frac{(a-1)g(t)^2 G(t)^{a-2}}{1-G(t)^a} + \frac{a g(t)^2 G(t)^{2a-2}}{[1-G(t)^a]^2}\right]$$

$$+ \frac{(m-1)ab[g'(t)G(t)^{a-1}[1-G(t)^a]^{b-1}]}{1-[1-G(t)^a]^b} + \frac{(m-1)ab[(a-1)g(t)^2 G(t)^{a-2}[1-G(t)^a]^{b-1}]}{1-[1-G(t)^a]^b}$$

$$- \frac{(m-1)ab[a(b-1)g(t)^2 G(t)^{2a-2}[1-G(t)^a]^{b-2}]}{1-[1-G(t)^a]^b}$$

$$- (m-1)\left[\frac{ab\,g(t)G(t)^{a-1}[1-G(t)^a]^{b-1}}{1-[1-G(t)^a]^b}\right]^2$$

The critical points of $h(t)$ are the roots of the equation

$$\frac{d\log[h(t;a,b,m,n)]}{dt}$$

$$= \frac{g'(t)}{g(t)} + (a-1)\frac{g(t)}{G(t)} + a(1-bn)\frac{g(t)G(t)^{a-1}}{1-G(t)^a} + (m-1)\frac{ab\,g(t)G(t)^{a-1}[1-G(t)^a]^{b-1}}{1-[1-G(t)^a]^b}$$

$$+ \frac{ab\,g(t)G(t)^{a-1}[1-G(t)^a]^{bn-1}[1-[1-G(t)^a]^b]^{m-1}}{B(m,n)[1 - I_{1-[1-G(t)^a]^b}(m,n)]} \qquad (19)$$

There may be more than one root to (19). If $t = t_0$ is a root of (19) then it corresponds to a local maximum, a local minimum or a point of inflexion depending on whether

$\omega(t_0) < 0$, $\omega(t_0) > 0$ or $\omega(t_0) = 0$ where $\omega(t) = \dfrac{d^2}{dt^2}\log[h(t;a,b,m,n)]$

$$\omega(t) = \frac{g(t)g''(t) - [g'(t)]^2}{g(t)^2} + (a-1)\frac{G(t)g'(t) - g(t)^2}{G(t)^2}$$

$$+ a(1-bn)\left[\frac{g'(t)G(t)^{a-1}}{1-G(t)^a} + \frac{(a-1)g(t)^2 G(t)^{a-2}}{1-G(t)^a} + \frac{a g(t)^2 G(t)^{2a-2}}{[1-G(t)^a]^2}\right]$$

$$+ \frac{(m-1)ab[g'(t)G(t)^{a-1}[1-G(t)^a]^{b-1}]}{1-[1-G(t)^a]^b} + \frac{(m-1)ab[(a-1)g(t)^2 G(t)^{a-2}[1-G(t)^a]^{b-1}]}{1-[1-G(t)^a]^b}$$

$$- \frac{(m-1)ab[a(b-1)g(t)^2 G(t)^{2a-2}[1-G(t)^a]^{b-2}]}{1-[1-G(t)^a]^b} - (m-1)\left[\frac{ab\,g(t)G(t)^{a-1}[1-G(t)^a]^{b-1}}{1-[1-G(t)^a]^b}\right]^2$$

$$+ \frac{ab\,g'(t)G(t)^{a-1}[1-G(t)^a]^{bn-1}[1-[1-G(t)^a]^b]^{m-1}}{B(m,n)\,[1 - I_{1-[1-G(t)^a]^b}(m,n)]}$$

$$+ \frac{ab(a-1)g(t)^2 G(t)^{a-2}[1-G(t)^a]^{bn-1}[1-[1-G(t)^a]^b]^{m-1}}{B(m,n)[1 - I_{1-[1-G(t)^a]^b}(m,n)]}$$

$$-\frac{a^2 b(bn-1) g(t)^2 G(t)^{2a-2}[1-G(t)^a]^{bn-2}[1-[1-G(t)^a]^b]^{m-1}}{B(m,n)[1-I_{1-[1-G(t)^a]^b}(m,n)]}$$

$$+\frac{a^2 b^2 (m-1) g(t)^2 G(t)^{2a-2}[1-G(t)^a]^{b(n+1)-2}[1-[1-G(t)^a]^b]^{m-2}}{B(m,n)[1-I_{1-[1-G(t)^a]^b}(m,n)]}$$

$$+\left\{\frac{ab\, g(t) G(t)^{a-1}[1-G(t)^a]^{bn-1}[1-[1-G(t)^a]^b]^{m-1}}{B(m,n)[1-I_{1-[1-G(t)^a]^b}(m,n)]}\right\}^2$$

## 5. Estimation

### 5.1 Maximum likelihood method

The model parameters of the $BKw-G$ distribution can be estimated by maximum likelihood. Let $t = (t_1, t_2, \ldots t_n)^T$ be a random sample of size $n$ from $BKw-G$ with parameter vector $\rho = (m, n, a, b, \boldsymbol{\beta}^T)^T$, where $\boldsymbol{\beta} = (\beta_1, \beta_2, \ldots \beta_q)^T$ corresponds to the parameter vector of the baseline distribution $G$. Then the log-likelihood function for $\rho$ is given by

$$\ell = \ell(\boldsymbol{\theta}) = r\log(ab) + \sum_{i=0}^{r}\log[g(t_i, \boldsymbol{\beta})] - r\log[B(m,n)] + (a-1)\sum_{i=0}^{r}\log[G(t_i, \boldsymbol{\beta})]$$

$$+ (bn-1)\sum_{i=0}^{r}\log[1-G(t_i, \boldsymbol{\beta})^a] + (m-1)\sum_{i=1}^{r}\log[1-[1-G(t_i, \boldsymbol{\beta})^a]^b] \qquad (20)$$

This log-likelihood function can not be solved analytically because of its complex form but it can be maximized numerically by employing global optimization methods available with software's like R, SAS, Mathematica or by solving the nonlinear likelihood equations obtained by differentiating (20).

By taking the partial derivatives of the log-likelihood function with respect to $m, n, a, b$ and $\boldsymbol{\beta}$ components of the score vector $U_{\boldsymbol{\theta}} = (U_m, U_n, U_a, U_b, U_{\boldsymbol{\beta}^T})^T$ can be obtained as follows:

$$U_m = \frac{\partial \ell}{\partial m} = -r\psi(m) + r\psi(m+n) + \sum_{i=1}^{r}\log[1-[1-G(t_i, \boldsymbol{\beta})^a]^b]$$

$$U_n = \frac{\partial \ell}{\partial n} = -r\psi(n) + r\psi(m+n) + b\sum_{i=0}^{r}\log[1-G(t_i, \boldsymbol{\beta})^a]$$

$$U_a = \frac{\partial \ell}{\partial a} = \frac{r}{a} + \sum_{i=0}^{r}\log[G(t_i, \boldsymbol{\beta})] + (1-bn)\sum_{i=0}^{r}\frac{G(t_i, \boldsymbol{\beta})^a \log[G(t_i, \boldsymbol{\beta})]}{1-G(t_i, \boldsymbol{\beta})^a}$$

$$+ (m-1)\sum_{i=1}^{r}\frac{b[1-G(t_i, \boldsymbol{\beta})^a]^{b-1} G(t_i, \boldsymbol{\beta})^a \log[G(t_i, \boldsymbol{\beta})]}{1-[1-G(t_i, \boldsymbol{\beta})^a]^b}$$

$$U_b = \frac{\partial \ell}{\partial b} = \frac{r}{b} + n\sum_{i=0}^{r}\log[1-G(t_i, \boldsymbol{\beta})^a] + (1-m)\sum_{i=0}^{r}\frac{[1-G(t_i, \boldsymbol{\beta})^a]^b \log[1-G(t_i, \boldsymbol{\beta})^a]}{1-[1-G(t_i, \boldsymbol{\beta})^a]^b}$$

$$U_\beta = \frac{\partial \ell}{\partial \boldsymbol{\beta}} = \sum_{i=0}^{r}\frac{g^{(\boldsymbol{\beta})}(t_i, \boldsymbol{\beta})}{g(t_i, \boldsymbol{\beta})} + (a-1)\sum_{i=0}^{r}\frac{G^{(\boldsymbol{\beta})}(t_i, \boldsymbol{\beta})}{G(t_i, \boldsymbol{\beta})} + (1-bn)\sum_{i=0}^{r}\frac{a G(t_i, \boldsymbol{\beta})^{a-1} G^{(\boldsymbol{\beta})}(t_i, \boldsymbol{\beta})}{1-G(t_i, \boldsymbol{\beta})^a}$$

$$+ (m-1) \sum_{i=0}^{r} \frac{b\,[1-G(t_i,\boldsymbol{\beta})^a]^{b-1}\,a\,G(t_i,\boldsymbol{\beta})^{a-1}\,G^{(\boldsymbol{\beta})}(t_i,\boldsymbol{\beta})}{1-[1-G(t_i,\boldsymbol{\beta})^a]^b}$$

Where $\psi(.)$ is the digamma function.

### 5.2 Asymptotic standard error and confidence interval for the mles

The asymptotic variance-covariance matrix of the MLEs of parameters can obtained by inverting the Fisher information matrix $I(\boldsymbol{\theta})$ which can be derived using the second partial derivatives of the log-likelihood function with respect to each parameter. The $ij^{th}$ elements of $I_n(\boldsymbol{\theta})$ are given by

$$I_{ij} = -E\left(\frac{\partial^2 l(\boldsymbol{\theta})}{\partial \theta_i \, \partial \theta_j}\right), \quad i,j = 1, 2, \cdots, 3+q$$

The exact evaluation of the above expectations may be cumbersome. In practice one can estimate $I_n(\boldsymbol{\theta})$ by the observed Fisher's information matrix $\hat{I}_n(\hat{\boldsymbol{\theta}})$ is defined as:

$$\hat{I}_{ij} \approx \left(-\frac{\partial^2 l(\boldsymbol{\theta})}{\partial \theta_i \, \partial \theta_j}\right)_{\boldsymbol{\theta}=\hat{\boldsymbol{\theta}}}, \quad i,j = 1, 2, \cdots, 3+q$$

Using the general theory of MLEs under some regularity conditions on the parameters as $n \to \infty$ the asymptotic distribution of $\sqrt{n}\,(\hat{\boldsymbol{\theta}} - \boldsymbol{\theta})$ is $N_k(0, V_n)$ where $V_n = (v_{jj}) = I_n^{-1}(\boldsymbol{\theta})$. The asymptotic behaviour remains valid if $V_n$ is replaced by $\hat{V}_n = \hat{I}^{-1}(\hat{\boldsymbol{\theta}})$. This This result can be used to provide large sample standard errors and also construct confidence intervals for the model parameters. Thus an approximate standard error and $(1-\gamma/2)100\%$ confidence interval for the mle of $j^{th}$ parameter $\theta_j$ are respectively given by $\sqrt{\hat{v}_{jj}}$ and $\hat{\theta}_j \pm Z_{\gamma/2}\sqrt{\hat{v}_{jj}}$, where $Z_{\gamma/2}$ is the $\gamma/2$ point of standard normal distribution.

As an illustration on the MLE method its large sample standard errors, confidence interval in the case of $BKw-E(m,n,a,b,\lambda)$ is discussed in an appendix.

### 5.3 Method of moments

Here an alternative method to estimation of the parameters is discussed. Since the moments are not in closed form, the estimation by the usual method of moments may not be tractable. Therefore following (Barreto-Souzai et al., 2013) we get

$$E\{[1-[1-G(t)^a]^b]\} = \int_{-\infty}^{\infty} \{[1-[1-G(t)^a]^b]\} \frac{a\,b\,g(t)G(t)^{a-1}[1-G(t)^a]^{bn-1}}{B(m,n)\,[1-[1-G(t)^a]^b]^{1-m}} dt$$

$$= \frac{1}{B(m,n)} \int_{-\infty}^{\infty} \frac{a\,b\,g(t)G(t)^{a-1}[1-G(t)^a]^{bn-1}}{[1-[1-G(t)^a]^b]^{-m}} dt$$

Let, $u = 1-[1-G(t)^a]^b$ then $du = ab\,g(t)G(t)^{a-1}[1-G(t)^a]^{b-1}\,dt$

$$= \frac{1}{B(m,n)}\int_0^1 u^m(1-u)^{n-1}\,du = \frac{1}{B(m,n)}\int_0^1 u^{(m+1)-1}(1-u)^{n-1}\,du = \frac{B(m+1,n)}{B(m,n)}$$

$$E\{[1-[1-G(t)^a]^b]^v\} = \int_{-\infty}^{\infty}\{[1-[1-G(t)^a]^b]^v\}\frac{ab\,g(t)G(t)^{a-1}[1-G(t)^a]^{bn-1}}{B(m,n)\,[1-[1-G(t)^a]^b]^{1-m}}\,dt$$

$$= \frac{1}{B(m,n)}\int_{-\infty}^{\infty}\frac{ab\,g(t)G(t)^{a-1}[1-G(t)^a]^{bn-1}}{[1-[1-G(t)^a]^b]^{1-m-v}}\,dt$$

Let, $u = 1-[1-G(t)^a]^b$ then $du = ab\,g(t)G(t)^{a-1}[1-G(t)^a]^{b-1}\,dt$

$$= \frac{1}{B(m,n)}\int_0^1 \frac{(1-u)^{n-1}}{u^{1-m-v}}\,du = \frac{1}{B(m,n)}\int_0^1 u^{v+m-1}(1-u)^{n-1}\,du = \frac{B(v+m,n)}{B(m,n)}$$

Therefore we have,

$$E[[1-[1-G(t)^a]^b]^v] = \frac{B(m+v,n)}{B(m,n)}, \text{ for } m,n>0,\ v=1,2,3,... \quad (21)$$

For a random sample $t_1, t_2, ..., t_n$ from a population with $BKw-G$ distribution, the model parameters can be estimated using (21) by numerically solving the equations

$$\frac{1}{n}\sum_{i=1}^{n}[[1-[1-G(t_i)^a]^b]^v] = \frac{B(m+v,n)}{B(m,n)}, \text{ for } m,n>0,\ v=1,2,3,...$$

## 6. Modelling applications

In this subsection, we consider real two data sets to establish the supremacy of the $BKw-G$ distributions over some of the distributions nested with in it when the Weibull distribution is chosen as $G$. We have used maximum likelihood method obtain the estimates of parameters. Additionally we have provided standard errors, and confidence interval (see Appendix) of the parameter estimators. Models are compared employing the Akaike Information Criterion (AIC). $AIC = 2k - 2l$ where $k$ is the number of parameters in the distribution and $l$ is the maximized value of the log-likelihood function under the considered model.

**Example I:**

The first data is about 346 nicotine measurements made from several brands of cigarettes in 1998. This data was originally collected by the Federal Trade Commission which is an independent agency of the US government, whose main mission is the promotion of consumer protection.[http://www.ftc.gov/reports/tobacco or http: //pw1.netcom.com/rdavis2/smoke. html.]. The data set is:

{1.3, 1.0, 1.2, 0.9, 1.1, 0.8, 0.5, 1.0, 0.7, 0.5, 1.7, 1.1, 0.8, 0.5, 1.2, 0.8, 1.1, 0.9, 1.2, 0.9, 0.8, 0.6, 0.3, 0.8, 0.6, 0.4, 1.1, 1.1, 0.2, 0.8, 0.5, 1.1, 0.1, 0.8, 1.7, 1.0, 0.8, 1.0, 0.8, 1.0, 0.2, 0.8, 0.4, 1.0, 0.2, 0.8, 1.4, 0.8, 0.5, 1.1, 0.9, 1.3, 0.9, 0.4, 1.4, 0.9, 0.5, 1.7, 0.9, 0.8, 0.8, 1.2, 0.9, 0.8, 0.5, 1.0, 0.6, 0.1, 0.2, 0.5, 0.1, 0.1, 0.9, 0.6, 0.9, 0.6, 1.2, 1.5, 1.1, 1.4, 1.2, 1.7, 1.4, 1.0, 0.7, 0.4, 0.9, 0.7, 0.8, 0.7, 0.4, 0.9, 0.6, 0.4, 1.2, 2.0, 0.7, 0.5, 0.9, 0.5, 0.9, 0.7, 0.9, 0.7, 0.4, 1.0, 0.7, 0.9, 0.7, 0.5, 1.3, 0.9, 0.8, 1.0, 0.7, 0.7, 0.6, 0.8, 1.1, 0.9, 0.9, 0.8, 0.8, 0.7, 0.7, 0.4, 0.5, 0.4, 0.9, 0.9, 0.7, 1.0, 1.0, 0.7, 1.3, 1.0, 1.1, 1.1, 0.9, 1.1, 0.8, 1.0, 0.7, 1.6, 0.8, 0.6, 0.8, 0.6, 1.2, 0.9, 0.6, 0.8, 1.0, 0.5, 0.8, 1.0, 1.1, 0.8, 0.8, 0.5, 1.1, 0.8, 0.9, 1.1, 0.8, 1.2, 1.1, 1.2, 1.1, 1.2, 0.2, 0.5, 0.7, 0.2, 0.5, 0.6, 0.1, 0.4, 0.6, 0.2, 0.5, 1.1, 0.8, 0.6, 1.1, 0.9, 0.6, 0.3, 0.9, 0.8, 0.8, 0.6, 0.4, 1.2, 1.3, 1.0, 0.6, 1.2, 0.9, 1.2, 0.9, 0.5, 0.8, 1.0, 0.7, 0.9, 1.0, 0.1, 0.2, 0.1, 0.1, 1.1, 1.0, 1.1, 0.7, 1.1, 0.7, 1.8, 1.2, 0.9, 1.7, 1.2, 1.3, 1.2, 0.9, 0.7, 0.7, 1.2, 1.0, 0.9, 1.6, 0.8, 0.8, 1.1, 1.1, 0.8, 0.6, 1.0, 0.8, 1.1, 0.8, 0.5, 1.5, 1.1, 0.8, 0.6, 1.1, 0.8, 1.1, 0.8, 1.5, 1.1, 0.8, 0.4, 1.0, 0.8, 1.4, 0.9, 0.9, 1.0, 0.9, 1.3, 0.8, 1.0, 0.5, 1.0, 0.7, 0.5, 1.4, 1.2, 0.9, 1.1, 0.9, 1.1, 1.0, 0.9, 1.2, 0.9, 1.2, 0.9, 0.5, 0.9, 0.7, 0.3, 1.0, 0.6, 1.0, 0.9, 1.0, 1.1, 0.8, 0.5, 1.1, 0.8, 1.2, 0.8, 0.5, 1.5, 1.5, 1.0, 0.8, 1.0, 0.5, 1.7, 0.3, 0.6, 0.6, 0.4, 0.5, 0.5, 0.7, 0.4, 0.5, 0.8, 0.5, 1.3, 0.9, 1.3, 0.9, 0.5, 1.2, 0.9, 1.1, 0.9, 0.5, 0.7, 0.5, 1.1, 1.1, 0.5, 0.8, 0.6, 1.2, 0.8, 0.4, 1.3, 0.8, 0.5, 1.2, 0.7, 0.5, 0.9, 1.3, 0.8, 1.2, 0.9}

**Table 1:** MLEs, standard errors and 95% confidence intervals (in parentheses) and the AIC values for the nicotine measurements data.

| Parameters | $B-W$ | $Kw-W$ | $BKw-W$ |
|---|---|---|---|
| $\hat{a}$ | --- | 0.792 (0.149) (0.49, 1.08) | 0.207 (0.282) (-0.35, 0.76) |
| $\hat{b}$ | --- | 0.430 (0.188) (0.06, 0.79) | 0.776 (0.191) (0.40, 1.15) |
| $\hat{m}$ | 0.777 (0.155) (0.47, 1.08) | --- | 2.647 (2.81) (-2.86, 8.15) |
| $\hat{n}$ | 2.027 (1.996) (-1.89, 5.94) | --- | 0.298 (0.900) (-1.47, 2.06) |
| $\hat{\lambda}$ | 0.431 (0.461) (-0.47, 1.33) | 2.326 (1.060) (0.25, 4.40) | 4.636 (0.534) (3.58, 5.68) |
| $\hat{\beta}$ | 3.158 (0.403) (2.37, 3.94) | 3.00 (0.259) (2.49, 3.51) | 2.912 (0.201) (2.52, 3.31) |
| log-likelihood($l_{\max}$) | -113.08 | -112.58 | **-110.06** |
| AIC | 234.16 | 233.16 | **232.12** |

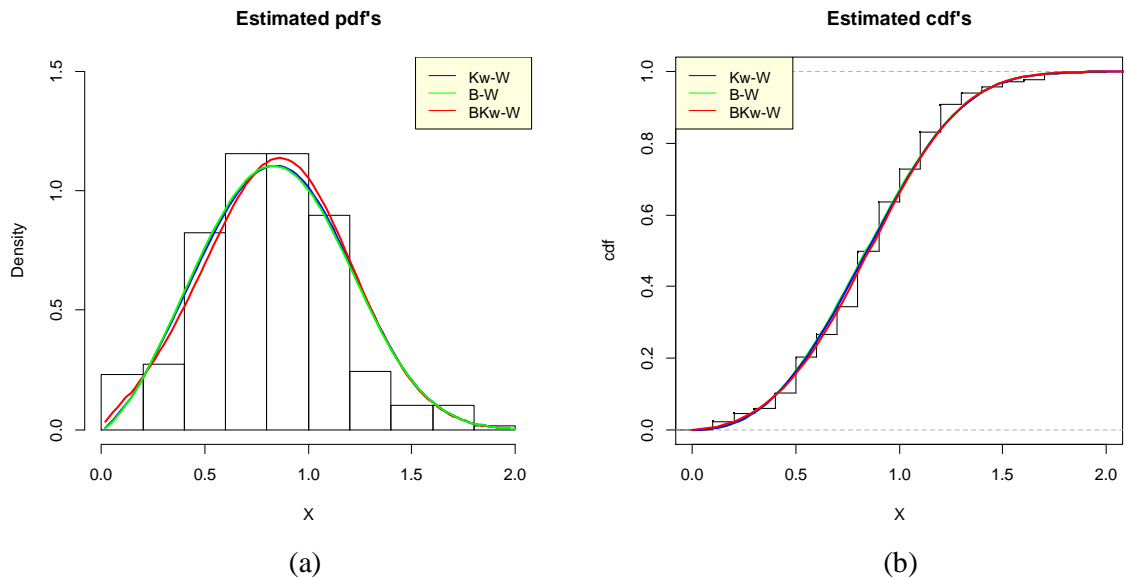

**Fig: 3** Plots of the (a) observed histogram and estimated pdf's and (b) estimated cdf's for the $Kw-W$, $B-W$ and $BKw-W$ for example I.

**Example II:**

The data set for the second application is a subset of the data reported by Bekker et al. (2000) which corresponds to the survival times (in years) of a group of patients given chemotherapy treatment alone. This data set consisting of survival times (in years) for 46 patients is:

{0.047, 0.115, 0.121, 0.132, 0.164, 0.197, 0.203, 0.260, 0.282, 0.296, 0.334, 0.395, 0.458, 0.466, 0.501, 0.507, 0.529, 0.534, 0.540, 0.641, 0.644, 0.696, 0.841, 0.863, 1.099, 1.219, 1.271, 1.326, 1.447, 1.485, 1.553, 1.581, 1.589, 2.178, 2.343, 2.416, 2.444, 2.825, 2.830, 3.578, 3.658, 3.743, 3.978, 4.003, 4.033}

**Table 2:** MLEs, standard errors and 95% confidence intervals (in parentheses) and the AIC values for survival times of cancer patients data.

| Parameters | $B-W$ | $Kw-W$ | $BKw-W$ |
|---|---|---|---|
| $\hat{a}$ | --- | 2.160 (1.584) (-0.94, 5.26) | 6.525 (0.079) (6.37, 6.67) |
| $\hat{b}$ | --- | 0.208 (0.204) (-0.19, 0.61) | 0.317 (0.323) (-0.32, 0.95) |
| $\hat{m}$ | 0.467 (0.807) (-1.11, 2.05) | --- | 5.183 (0.004) (5.17, 5.19) |
| $\hat{n}$ | 0.695 (0.555) (-0.39, 1.78) | --- | 1.388 (0.002) (1.38, 1.39) |
| $\hat{\lambda}$ | 2.024 (2.856) (-3.57, 7.62) | 4.521 (3.900) (-3.12, 12.16) | 0.223 (0.054) (0.12, 0.33) |
| $\hat{\beta}$ | 3.418 (6.538) (-9.39, 16.2) | 0.836 (0.196) (0.45, 1.22) | 0.224 (0.189) (-0.15, 0.59) |
| log-likelihood($l_{\max}$) | -58.07 | -57.72 | **-55.46** |
| AIC | 124.14 | 123.44 | **122.92** |

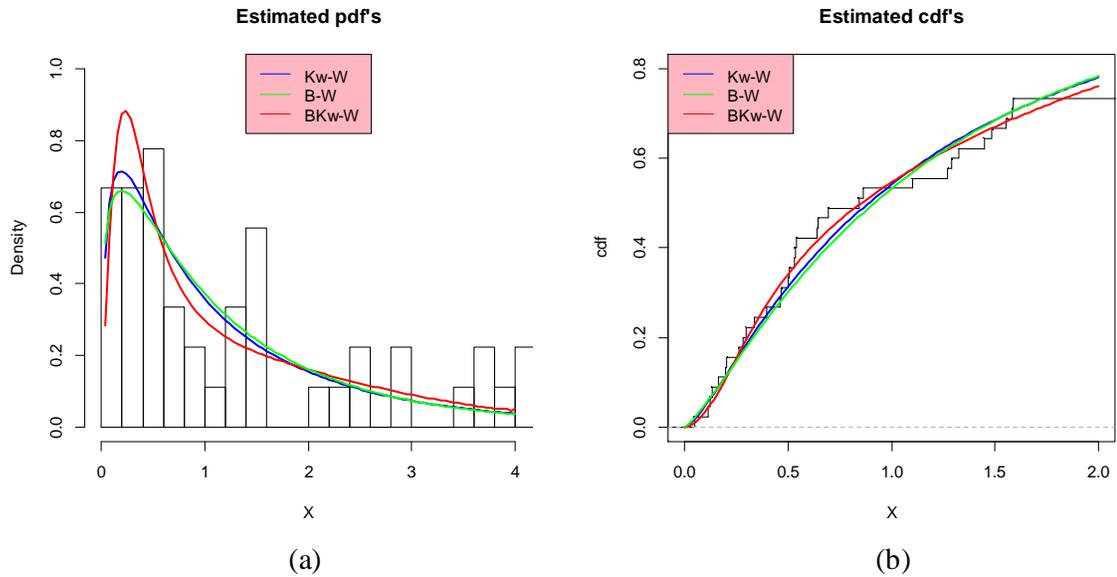

(a)      (b)

**Fig: 4** Plots of the (a) observed histogram and estimated pdf's and (b) estimated cdf's for the $Kw-W$, $B-W$ and $BKw-W$ for example II.

In the Table 1(Table 2), the estimates, MLEs, standard errors (SEs) (in parentheses) and confidence intervals (in parentheses) of the parameters for all the fitted models along with their AIC values for example I (example II) are presented. It is evident that in both the examples the $BKw-W$ distribution is a better model than both its sub models $Kw-W$ and $B-W$ with respect to values of the AIC. This is further validated by providing plots of fitted densities with histogram of the observed data in figure 3 (a) (figure 4 (a)) and fitted cdfs with cumulative frequency of observed data in figure 3 (b) (figure 4 (b)) for example I (example II). These plots indicate that the proposed distributions provide closer fit to both the observed data sets.

## 7. Some other classes of Beta generated distributions

Proceeding in the similar manner we propose to investigate the following classes distributions which will be reported in due time.

I. Beta Marshall-Olkin Kumaraswamy-$G$ family of distributions ($BMOKw-G$) with pdf

$$f^{BMOKwG}(t;m,n,l,a,b) = \frac{l^n a b\, g(t) G(t)^{a-1}[1-G(t)^a]^{bn-1}[1-[1-G(t)^a]^b]^{m-1}}{B(m,n)[1-(1-l)[1-G(t)^a]^b]^{m+n}}$$

$$0<t<\infty,\ 0<a,b<\infty,\ m,n>0,\ l>0$$

II. Beta Kumaraswamy Marshall-Olkin-$G$ family of distributions ($BKwMO-G$) with pdf

$$f^{BKwMOG}(t;a,b,m,n,\alpha) = \frac{ab\alpha\, g(t)\, G(t)^{a-1}}{B(m,n)[1-\overline{\alpha}\,\overline{G}(t)]^{a+1}}\,[\,1-[G(t)/\{1-\overline{\alpha}\,\overline{G}(t)\}]^a\,]^{bn-1}$$

$$\times [1-[1-[G(t)/\{1-\overline{\alpha}\,\overline{G}(t)\}]^a]^b]^{m-1}$$

$$0<t<\infty,\ 0<a,b<\infty,\ m,n>0,\alpha>0$$

III. Beta Kumaraswamy Generalized Marshall-Olkin-$G$ family of distributions ($BKwGMO-G$) with pdf:

$$f^{BKwGMOG}(t;a,b,m,n,\alpha,\theta)$$

$$= \frac{1}{B(m,n)} \frac{ab\theta\alpha^\theta\, g(t)\overline{G}(t)^{\theta-1}}{[1-\overline{\alpha}\,\overline{G}(t)]^{\theta+1}} \left\{ 1-\left[\frac{\alpha\,\overline{G}(t)}{1-\overline{\alpha}\,\overline{G}(t)}\right]^\theta \right\}^{a-1}$$

$$\times \left\{ 1-\left\{1-\left[\frac{\alpha\,\overline{G}(t)}{1-\overline{\alpha}\,\overline{G}(t)}\right]^\theta\right\}^a \right\}^{bn-1} \left[ 1-\left\{1-\left\{1-\left[\frac{\alpha\,\overline{G}(t)}{1-\overline{\alpha}\,\overline{G}(t)}\right]^\theta\right\}^a\right\}^b \right]^{m-1}$$

$$0<t<\infty,\ 0<a,b<\infty,\ m,n>0,\ \alpha>0,\ \theta>0$$

IV. Beta Generalized Marshall-Olkin Kumaraswamy-$G$ family of distributions $BGMOKw-G$ with pdf

$$f^{BGMOKwG}(t;m,n,a,b,\alpha,\theta) = \frac{\theta \alpha^\theta \, a b \, g(t) G(t)^{a-1} [1-G(t)^a]^{b-1} [[1-G(t)^a]^b]^{\theta-1}}{B(m,n)[1-\bar{\alpha}[1-G(t)^a]^b]^{\theta+1}}$$

$$\times \left[1 - \left[\frac{\alpha[1-G(t)^a]^b}{1-\bar{\alpha}[1-G(t)^a]^b}\right]^\theta\right]^{m-1} \left[\left[\frac{\alpha[1-G(t)^a]^b}{1-\bar{\alpha}[1-G(t)^a]^b}\right]^\theta\right]^{n-1}$$

$$0 < t < \infty, \ 0 < a,b < \infty, \ m,n > 0, \alpha > 0, \theta > 0$$

## 8. Concluding discussion and remarks

A new beta extended Kumaraswamy generalized family of distributions which includes three recently proposed new families of distributions namely the Garhy generated family (*Garhy et al.*, 2016), Beta-Dagum distribution and Beta-Singh-Maddala distribution (Domma and Condino, 2016) is introduced and some of its important properties are studied. The maximum likelihood and moment method for estimating the parameters are also discussed. Comparative data modeling application of the proposed model with some of its sub models is carried out considering two data sets reveals its superiority. Some more beta generated families are introduced for future investigation.

### Appendix: Maximum likelihood method for $BKw-E$ distribution:

The $BKw-E$ distribution is a special case of $BKw-G$ with pdf

$$f^{BKwE}(t;m,n,a,b,\lambda) = \frac{1}{B(m,n)} a b \lambda e^{-\lambda t} \{1-e^{-\lambda t}\}^{a-1} [1-\{1-e^{-\lambda t}\}^a]^{bn-1} [1-[1-\{1-e^{-\lambda t}\}^a]^b]^{m-1}$$

for $m > 0, n > 0, \lambda > 0, \ a > 0, b > 0, t > 0$

For a random sample of size *n*, the log-likelihood function for the parameter vector $\boldsymbol{\theta} = (m,n,a,b,\lambda)^T$ is given by

$$\ell = \ell(\boldsymbol{\theta}) = r \log(ab) + r \log \lambda - \lambda \sum_{i=0}^{r} t_i - r \log[B(m,n)] + (a-1) \sum_{i=0}^{r} \log(1-e^{-\lambda t_i})$$

$$+ (bn-1) \sum_{i=0}^{r} \log[1-(1-e^{-\lambda t_i})^a] + (m-1) \sum_{i=0}^{r} \log[1-\{1-(1-e^{-\lambda t_i})^a\}^b]$$

The corresponding components of the score vector $\boldsymbol{\theta} = (m,n,a,b,\lambda)^T$ are

$$U_m = \frac{\partial \ell}{\partial m} = -r\psi(m) - r\psi(m+n) + \sum_{i=0}^{r} \log[1-\{1-(1-e^{-\lambda t_i})^a\}^b]$$

$$U_n = \frac{\partial \ell}{\partial n} = -r\psi(n) - r\psi(m+n) + b \sum_{i=0}^{r} \log[1-(1-e^{-\lambda t_i})^a]$$

$$U_a = \frac{\partial \ell}{\partial a} = \frac{r}{a} + \sum_{i=0}^{r} \log(1-e^{-\lambda t_i}) + (1-bn)\sum_{i=0}^{r} \frac{(1-e^{-\lambda t_i})^a \log(1-e^{-\lambda t_i})}{1-(1-e^{-\lambda t_i})^a}$$

$$+ (m-1)\sum_{i=0}^{r} \frac{b[1-(1-e^{-\lambda t_i})^a]^{b-1}(1-e^{-\lambda t_i})^a \log[(1-e^{-\lambda t_i})]}{1-[1-(1-e^{-\lambda t_i})^a]^b}$$

$$U_b = \frac{\partial \ell}{\partial b} = \frac{r}{b} + n\sum_{i=0}^{r}\log[1-(1-e^{-\lambda t_i})^a] + (1-m)\sum_{i=0}^{r}\frac{[1-(1-e^{-\lambda t_i})^a]^b \log[1-(1-e^{-\lambda t_i})^a]}{1-[1-(1-e^{-\lambda t_i})^a]^b}$$

$$U_\lambda = \frac{\partial \ell}{\partial \lambda} = \frac{r}{\lambda} - \sum_{i=0}^{r}t_i + (a-1)\sum_{i=0}^{r}\frac{\lambda e^{-\lambda t_i}}{1-e^{-\lambda t_i}} + (1-bn)\sum_{i=0}^{r}\frac{a\lambda(1-e^{-\lambda t_i})^{a-1}e^{-\lambda t_i}}{1-(1-e^{-\lambda t_i})^a}$$

$$+ (m-1)\sum_{i=0}^{n} \frac{ab\lambda[1-(1-e^{-\lambda t_i})^a]^{b-1}(1-e^{-\lambda t_i})^{a-1}e^{-\lambda t_i}}{1-[1-(1-e^{-\lambda t_i})^a]^b}$$

The asymptotic variance covariance matrix for mles of the parameters of $BKw-E$ $(m,n,a,b,\lambda)$ distribution is estimated by

$$\hat{I}_n^{-1}(\hat{\theta}) = \begin{pmatrix} \text{var}(\hat{m}) & \text{cov}(\hat{m},\hat{n}) & \text{cov}(\hat{m},\hat{a}) & \text{cov}(\hat{m},\hat{b}) & \text{cov}(\hat{m},\hat{\lambda}) \\ \text{cov}(\hat{n},\hat{m}) & \text{var}(\hat{n}) & \text{cov}(\hat{n},\hat{a}) & \text{cov}(\hat{n},\hat{b}) & \text{cov}(\hat{n},\hat{\lambda}) \\ \text{cov}(\hat{a},\hat{m}) & \text{cov}(\hat{a},\hat{n}) & \text{var}(\hat{a}) & \text{cov}(\hat{a},\hat{b}) & \text{cov}(\hat{a},\hat{\lambda}) \\ \text{cov}(\hat{b},\hat{m}) & \text{cov}(\hat{b},\hat{n}) & \text{cov}(\hat{b},\hat{a}) & \text{var}(\hat{b}) & \text{cov}(\hat{b},\hat{\lambda}) \\ \text{cov}(\hat{\lambda},\hat{m}) & \text{cov}(\hat{\lambda},\hat{n}) & \text{cov}(\hat{\lambda},\hat{a}) & \text{cov}(\hat{\lambda},\hat{b}) & \text{var}(\hat{\lambda}) \end{pmatrix}$$

Where the elements of the information matrix $\hat{I}_n(\hat{\theta}) = \left(-\frac{\partial^2 l(\theta)}{\partial \theta_i \partial \theta_j}\right)_{\theta=\hat{\theta}}$

can be derived using the following second partial derivatives:

$$\frac{\partial^2 \ell}{\partial m^2} = -r\psi'(m) - r\psi'(m+n)$$

$$\frac{\partial^2 \ell}{\partial n^2} = -r\psi'(n) - r\psi'(m+n)$$

$$\frac{\partial^2 \ell}{\partial a^2} = -\frac{r}{a^2} + (1-bn)\sum_{i=0}^{r}\left(\frac{(1-e^{-\lambda t_i})^{2a}\log(1-e^{-\lambda t_i})^2}{\{1-(1-e^{-\lambda t_i})^a\}^2} + \frac{(1-e^{-\lambda t_i})^a \log(1-e^{-\lambda t_i})^2}{1-(1-e^{-\lambda t_i})^a}\right)$$

$$+ (1-m)\sum_{i=0}^{r}\frac{b^2[1-(1-e^{-\lambda t_i})^a]^{2(b-1)}(1-e^{-\lambda t_i})^{2a}\log[(1-e^{-\lambda t_i})]^2}{[1-[1-(1-e^{-\lambda t_i})^a]^b]^2}$$

$$+ (1-m)\sum_{i=0}^{r}\frac{b(b-1)[1-(1-e^{-\lambda t_i})^a]^{b-2}(1-e^{-\lambda t_i})^{2a}\log[(1-e^{-\lambda t_i})]^2}{1-[1-(1-e^{-\lambda t_i})^a]^b}$$

$$- (1-m)\sum_{i=0}^{r}\frac{b[1-(1-e^{-\lambda t_i})^a]^{b-1}(1-e^{-\lambda t_i})^a \log[(1-e^{-\lambda t_i})]^2}{1-[1-(1-e^{-\lambda t_i})^a]^b}$$

$$\frac{\partial^2 \ell}{\partial b^2} = \frac{r}{b} + (1-m)\sum_{i=0}^{r} \frac{[1-(1-e^{-\lambda t_i})^a]^{2b} \log[1-(1-e^{-\lambda t_i})^a]^2}{[1-[1-(1-e^{-\lambda t_i})^a]^b]^2}$$

$$+ (1-m)\sum_{i=0}^{r} \frac{[1-(1-e^{-\lambda t_i})^a]^b \log[1-(1-e^{-\lambda t_i})^a]^2}{1-[1-(1-e^{-\lambda t_i})^a]^b}$$

$$\frac{\partial^2 \ell}{\partial \lambda^2} = -\frac{r}{\lambda^2} + (a-1)\sum_{i=0}^{r}\left(\frac{e^{-2\lambda t_i} t_i^2}{(1-e^{-\lambda t_i})^2} - \frac{e^{-\lambda t_i} t_i^2}{1-e^{-\lambda t_i}}\right) + (1-bn)\sum_{i=0}^{r} \frac{a^2 (1-e^{-\lambda t_i})^{2(a-1)} e^{-2\lambda t_i} t_i^2}{\{1-(1-e^{-\lambda t_i})^a\}^2}$$

$$+ (1-bn)\sum_{i=0}^{r} \frac{a(a-1)(1-e^{-\lambda t_i})^{a-2} e^{-2\lambda t_i} t_i^2}{1-(1-e^{-\lambda t_i})^a} - (1-bn)\sum_{i=0}^{r} \frac{a(1-e^{-\lambda t_i})^{a-1} e^{-\lambda t_i} t_i^2}{1-(1-e^{-\lambda t_i})^a}$$

$$+ (1-m)\sum_{i=0}^{n} \frac{a^2 b^2 [1-(1-e^{-\lambda t_i})^a]^{2(b-1)} (1-e^{-\lambda t_i})^{2(a-1)} e^{-2\lambda t_i} t_i^2}{[1-[1-(1-e^{-\lambda t_i})^a]^b]^2}$$

$$+ (1-m)\sum_{i=0}^{n} \frac{a^2 b(b-1) [1-(1-e^{-\lambda t_i})^a]^{b-2} (1-e^{-\lambda t_i})^{2(a-1)} e^{-2\lambda t_i} t_i^2}{1-[1-(1-e^{-\lambda t_i})^a]^b}$$

$$- (1-m)\sum_{i=0}^{n} \frac{a(a-1)b [1-(1-e^{-\lambda t_i})^a]^{b-1} (1-e^{-\lambda t_i})^{a-2} e^{-2\lambda t_i} t_i^2}{1-[1-(1-e^{-\lambda t_i})^a]^b}$$

$$+ (1-m)\sum_{i=0}^{n} \frac{ab [1-(1-e^{-\lambda t_i})^a]^{b-1} (1-e^{-\lambda t_i})^{a-1} e^{-\lambda t_i} t_i^2}{1-[1-(1-e^{-\lambda t_i})^a]^b}$$

$$\frac{\partial^2 \ell}{\partial m \partial n} = r\psi'(m+n)$$

$$\frac{\partial^2 \ell}{\partial m \partial a} = \sum_{i=0}^{r} \frac{b[1-(1-e^{-\lambda t_i})^a]^{b-1} (1-e^{-\lambda t_i})^a \log[(1-e^{-\lambda t_i})]}{1-[1-(1-e^{-\lambda t_i})^a]^b}$$

$$\frac{\partial^2 \ell}{\partial m \partial b} = \sum_{i=0}^{r} -\frac{[1-(1-e^{-\lambda t_i})^a]^b \log[1-(1-e^{-\lambda t_i})^a]}{1-[1-(1-e^{-\lambda t_i})^a]^b}$$

$$\frac{\partial^2 \ell}{\partial m \partial \lambda} = \sum_{i=0}^{n} \frac{ab [1-(1-e^{-\lambda t_i})^a]^{b-1} (1-e^{-\lambda t_i})^{a-1} e^{-\lambda t_i} t_i}{1-[1-(1-e^{-\lambda t_i})^a]^b}$$

$$\frac{\partial^2 \ell}{\partial n \partial a} = b\sum_{i=0}^{r} -\frac{(1-e^{-\lambda t_i})^a \log(1-e^{-\lambda t_i})}{1-(1-e^{-\lambda t_i})^a}$$

$$\frac{\partial^2 \ell}{\partial n \partial m} = r\psi'(m+n)$$

$$\frac{\partial^2 \ell}{\partial n \partial b} = \sum_{i=0}^{r} \log[1-(1-e^{-\lambda t_i})^a]; \quad \frac{\partial^2 \ell}{\partial n \partial \lambda} = b\sum_{i=0}^{r} -\frac{a\lambda (1-e^{-\lambda t_i})^{a-1} e^{-\lambda t_i}}{1-(1-e^{-\lambda t_i})^a}$$

$$\frac{\partial^2 \ell}{\partial a \partial b} = -n\sum_{i=0}^{r} \frac{(1-e^{-\lambda t_i})^a \log(1-e^{-\lambda t_i})}{1-(1-e^{-\lambda t_i})} + (m-1)\sum_{i=0}^{r} \frac{[1-(1-e^{-\lambda t_i})^a]^{b-1} (1-e^{-\lambda t_i})^a \log(1-e^{-\lambda t_i})}{1-[1-(1-e^{-\lambda t_i})^a]^b}$$

$$+ b(m-1) \sum_{i=0}^{r} \frac{[1-(1-e^{-\lambda t_i})^a]^{2b-1} (1-e^{-\lambda t_i})^a \log(1-e^{-\lambda t_i}) \log[1-(1-e^{-\lambda t_i})^a]}{[1-[1-(1-e^{-\lambda t_i})^a]^b]^2}$$

$$+ b(m-1) \sum_{i=0}^{r} \frac{[1-(1-e^{-\lambda t_i})^a]^{b-1} (1-e^{-\lambda t_i})^a \log(1-e^{-\lambda t_i}) \log[1-(1-e^{-\lambda t_i})^a]}{1-[1-(1-e^{-\lambda t_i})^a]^b}$$

$$\frac{\partial^2 \ell}{\partial a \partial \lambda} = \sum_{i=0}^{r} \frac{\lambda e^{-\lambda t_i} t_i}{1-e^{-\lambda t_i}} + (1-bn) \sum_{i=0}^{r} \frac{a(1-e^{-\lambda t_i})^{a-1} e^{-\lambda t_i} t_i}{1-(1-e^{-\lambda t_i})^a}$$

$$+ (1-bn) \sum_{i=0}^{r} \frac{a(1-e^{-\lambda t_i})^{2a-1} e^{-\lambda t_i} \log(1-e^{-\lambda t_i}) t_i}{[1-(1-e^{-\lambda t_i})^a]^2} + (1-bn) \sum_{i=0}^{r} \frac{a(1-e^{-\lambda t_i})^{a-1} e^{-\lambda t_i} \log(1-e^{-\lambda t_i}) t_i}{1-(1-e^{-\lambda t_i})^a}$$

$$+ (m-1) \sum_{i=0}^{n} \frac{ab[1-(1-e^{-\lambda t_i})^a]^{b-1} (1-e^{-\lambda t_i})^{a-1} e^{-\lambda t_i} t_i}{1-[1-(1-e^{-\lambda t_i})^a]^b}$$

$$+ (m-1) \sum_{i=0}^{n} \frac{ab^2 [1-(1-e^{-\lambda t_i})^a]^{2(b-1)} (1-e^{-\lambda t_i})^{2a-1} e^{-\lambda t_i} t_i \log(1-e^{-\lambda t_i})}{[1-[1-(1-e^{-\lambda t_i})^a]^b]^2}$$

$$+ (m-1) \sum_{i=0}^{n} \frac{ab(b-1)[1-(1-e^{-\lambda t_i})^a]^{b-2} (1-e^{-\lambda t_i})^{2a-1} e^{-\lambda t_i} t_i \log(1-e^{-\lambda t_i})}{1-[1-(1-e^{-\lambda t_i})^a]^b}$$

$$- (m-1) \sum_{i=0}^{n} \frac{ab[1-(1-e^{-\lambda t_i})^a]^{b-1} (1-e^{-\lambda t_i})^{a-1} e^{-\lambda t_i} t_i \log(1-e^{-\lambda t_i})}{1-[1-(1-e^{-\lambda t_i})^a]^b}$$

$$\frac{\partial^2 \ell}{\partial b \partial \lambda} = -n \sum_{i=0}^{r} \frac{a(1-e^{-\lambda t_i})^{a-1} e^{-\lambda t_i} t_i}{1-(1-e^{-\lambda t_i})^a}$$

$$\frac{\partial \ell}{\partial \lambda} = \frac{r}{\lambda} - \sum_{i=0}^{r} t_i + (a-1) \sum_{i=0}^{r} \frac{\lambda e^{-\lambda t_i}}{1-e^{-\lambda t_i}} + (1-bn) \sum_{i=0}^{r} \frac{a\lambda(1-e^{-\lambda t_i})^{a-1} e^{-\lambda t_i}}{1-(1-e^{-\lambda t_i})^a}$$

$$+ (m-1) \sum_{i=0}^{n} \frac{ab[1-(1-e^{-\lambda t_i})^a]^{b-1} (1-e^{-\lambda t_i})^{a-1} e^{-\lambda t_i} t_i}{1-[1-(1-e^{-\lambda t_i})^a]^b}$$

$$+ (m-1) \sum_{i=0}^{n} \frac{ab[1-(1-e^{-\lambda t_i})^a]^{2b-1} (1-e^{-\lambda t_i})^{a-1} e^{-\lambda t_i} t_i \log[1-(1-e^{-\lambda t_i})^a]}{[1-[1-(1-e^{-\lambda t_i})^a]^b]^2}$$

$$+ (m-1) \sum_{i=0}^{n} \frac{ab[1-(1-e^{-\lambda t_i})^a]^{b-1} (1-e^{-\lambda t_i})^{a-1} e^{-\lambda t_i} t_i \log[1-(1-e^{-\lambda t_i})^a]}{1-[1-(1-e^{-\lambda t_i})^a]^b}$$

Where $\psi'(.)$ is the derivative of the digamma function.